\documentclass[graybox]{SNmult}
\usepackage{type1cm}        
\usepackage{makeidx}         
\usepackage{graphicx}        
\usepackage{multicol}        
\usepackage[bottom]{footmisc}

\usepackage{cite}
\usepackage{bm}
\usepackage{newtxtext}       %
\usepackage[varvw]{newtxmath}       
\usepackage{booktabs,multirow} 

\makeindex             
\newcommand\eref[1]{(\ref{#1})}
\newcommand*\xbar[1]{%
  \hbox{%
    \vbox{%
      \hrule height 0.5pt 
      \kern0.4ex
      \hbox{%
        \kern-0.05em
        \ensuremath{#1}%
        \kern-0.00em
      }%
    }%
  }%
}

\titlerunning{New Scheme Adaption Strategy}

\usepackage{amsmath}
\makeatletter
\AtBeginDocument{%
\@addtoreset{equation}{chapter}%
}
\makeatother

\newcommand{\bmF}{\bm{\mathcal{F}}}
\newcommand{\bmG}{\bm{\mathcal{G}}}
\newcommand{\mF}{\bm{F}}

\newcommand{\mG}{\bm{G}}

\newcommand{\mU}{\bm{U}}

\newcommand{\dx}{\Delta x}
\newcommand{\dy}{\Delta y}

\newcommand{\hf}{{\frac{1}{2}}}

\newcommand{\jph}{{j+\frac{1}{2}}}
\newcommand{\jmh}{{j-\frac{1}{2}}}
\newcommand{\kph}{{k+\frac{1}{2}}}
\newcommand{\kmh}{{k-\frac{1}{2}}}

\def\softd{{\leavevmode\setbox1=\hbox{d}%
          \hbox to 1.05\wd1{d\kern-0.4ex{\char039}\hss}}}

\begin{document}
\title*{New Scheme Adaption Strategy for Hyperbolic Conservation Laws}
\author{Shaoshuai Chu,~ Michael Herty,~ and~ Alexander Kurganov}
\institute{Shaoshuai Chu\at Department of Mathematics, RWTH Aachen University, 52056 Aachen, Germany, \email{chu@igpm.rwth-aachen.de}\and
Michael Herty\at Department of Mathematics, RWTH Aachen University, 52056 Aachen, Germany; Department of Mathematics and Applied
Mathematics, University of Pretoria, Private Bag X20, Hatfield 0028, South Africa, \email{herty@igpm.rwth-aachen.de}\and Alexander Kurganov
\at Department of Mathematics and Shenzhen International Center for Mathematics, Southern University of Science and Technology, Shenzhen,
518055, China, \email{alexander@sustech.edu.cn}}

\maketitle


\abstract{~We introduce a new scheme adaption strategy for one- and two-dimensional hyperbolic systems of conservation laws. The proposed
approach builds upon the adaptive framework introduced in [{\sc S. Chu, A. Kurganov, and I. Menshov}, Appl. Numer. Math., 209 (2025), pp.
155--170], where we first employed the smoothness indicator from [{\sc R. L\"ohner}, Comput. Methods. Appl. Mech. Eng., 61 (1987), pp.
323--338] to automatically detect ``rough'' and smooth parts of the computed solution, and then used different limiters in the detected
regions. This adaptive strategy was based on a threshold needed to sharply separate ``rough'' and smooth regions.\\[1.ex] 
In this paper, we propose a different adaption strategy. We use SBM-type limiters and vary one of the limiting parameters continuously to
allow a smooth transition between the ``rough'' and smooth areas. This way, compressive and overcompressive limiters are activated in the
shock and contact wave vicinities only, while we gradually switch to dissipative limiters in the smooth regions. A series of one- and
two-dimensional numerical tests for the Euler equations of gas dynamics demonstrates that the new scheme adaption strategy leads to a higher
resolution and reduced numerical dissipation.
}

\section{Introduction}
This paper focuses on the development of a new scheme adaption strategy for hyperbolic systems of conservation laws, which, in the
one-dimensional (1-D) and two-dimensional (2-D) cases, read as
\begin{equation}
\mU_t+\mF(\mU)_x=\bm0,
\label{1.1}
\end{equation}
and
\begin{equation}
\mU_t+\mF(\mU)_x+\mG(\mU)_y=\bm0,
\label{1.2}
\end{equation}
respectively. Here, $x$ and $y$ are spatial variables, $t$ is time, $\mU\in\mathbb R^d$ is the vector of unknowns, and
$\mF,\mG:\mathbb R^d\to\mathbb R^d$ are nonlinear flux functions.

It is well-known that, even for smooth initial data, solutions of nonlinear hyperbolic systems may develop complex structures, including
shocks, rarefactions, and contact discontinuities. This makes the construction of accurate and robust shock-capturing methods a challenging
task. In this work, we restrict our consideration to adaptive methods, which automatically detect ``rough'' parts of the computed solution
and apply different numerical procedures there and in the smooth regions. In particular, we build upon the adaptive strategy introduced in
\cite{CKM2025}, where the low-dissipation central-upwind (LDCU) numerical fluxes from \cite{CKX_24} were combined with a piecewise linear
reconstruction procedure. The slopes of the linear pieces were computed adaptively. To this end, ``rough'' regions were detected with the
help of the smoothness indicator (SI) proposed in \cite{RL87}. After that, the two-parameter family of the SBM-type limiters from
\cite{Lie03} were used with the parameters corresponding to either an overcompressive limiter (in the ``rough'' areas) or a dissipative one
(in the ``smooth'' areas). The latter was necessary to suppress artificial kinks, staircase-like structures, and spurious jump
discontinuities that may contaminate smooth solution profiles if they are captured using overcompressive limiters; see \cite{Lie03}. Since
the limiters utilized in \cite{CKM2025} were quite sharp (even the dissipative Minmod2 limiter), they were implemented in local
characteristic variables using a local characteristic decomposition.  

The main objective of the present paper is to improve the adaptive schemes introduced in \cite{CKM2025} by modifying the scheme adaption
procedure. In the conventional approach (also used in \cite{CKM2025}), the ``rough'' areas are identified as a group of cells (or cell
interfaces), in which the values of the SI are above a prescribed threshold. After that, a decision on which sets of parameters in the
SBM-type limiter should be used in ``rough'' and ``smooth'' areas is made: this completes the design of a scheme adaption method. 

An obvious drawback of such an approach is the sharp separation between the ``rough'' and ``smooth'' regions, which often leads to undesired
numerical artifacts. If the threshold is chosen too small, the overcompressive limiter is activated in excessively large parts of the
computational domain, which may result in artificial kinks, staircase-like structures, and spurious jump discontinuities appearing in the
areas where the exact solution is smooth. On the other hand, if the threshold is chosen too large, the dissipative limiter is used in
excessively large parts of the computational domain, and captured discontinuities, especially contact ones, may be oversmeared. To overcome
these drawbacks, we introduce transition layers between the ``rough'' and ``smooth'' areas by smoothly switching the parameters in the
SBM-type limiters. This allows us to keep using very sharp, overcompressive limiters right at the discontinuities and their immediate
neighborhood, while switching to compressive and dissipative limiters in the transition layers (this helps to keep the discontinuities much
better resolved without risking smooth solution structures to be affected), and then using dissipative limiters inside the ``smooth''
regions.  

The rest of the paper is organized as follows. In \S\ref{sec2}, we briefly review the 1-D LDCU scheme from \cite{CKX_24} and then introduce
the new scheme adaption strategy. In \S\ref{sec3}, we extend the proposed scheme adaption strategy to the 2-D case. In \S\ref{sec4}, we
apply the developed methods to a series of 1-D and 2-D test problems for the Euler equations of gas dynamics and demonstrate that the new
scheme adaption strategies lead to higher resolution compared with the resolution achieved by their counterparts from \cite{CKM2025}.

\section{One-Dimensional Scheme Adaption Algorithm}\label{sec2}
In this section, we consider the 1-D conservation law \eref{1.1} and describe the proposed adaptive algorithm. We will first briefly review
the 1-D LDCU scheme from \cite{CKX_24} as well as the SBM-type limiters (\S\ref{sec21}) and then introduce the new 1-D adaptive scheme
(\S\ref{sec22}).

\subsection{1-D Low-Dissipation Central-Upwind (LDCU) Scheme}\label{sec21}
Assume that the computational domain is partitioned into uniform cells $I_j:=[x_\jmh,x_\jph]$ of size $x_\jph-x_\jmh\equiv\dx$, centered at
$x_j:=\big(x_\jmh+x_\jph\big)/2$. We denote by $\,\xbar\mU_j(t)$ the cell averages of $\mU(\cdot,t)$ over $I_j$, namely,
\begin{equation*}
\xbar\mU_j(t)\approx\frac{1}{\dx}\int\limits_{I_j}\mU(x,t)\,{\rm d}x,
\end{equation*}
and assume that, at a given time level $t\ge0$, the values $\big\{\,\xbar\mU_j(t)\big\}$ are available. Note that most of the indexed
quantities depend on $t$, but this dependence will be suppressed hereafter for the sake of brevity.

In the semi-discrete LDCU scheme, the cell averages are evolved by solving the ODE system
\begin{equation}
\frac{{\rm d}\,\xbar\mU_j}{{\rm d}t}=-\frac{\bm{{\cal F}}_\jph-\bm{{\cal F}}_\jmh}{\dx},
\label{2.1}
\end{equation}
where $\bm{{\cal F}}_\jph=\bmF\big(\mU^-_\jph,\mU^+_\jph\big)$ is the LDCU numerical flux. In \cite{CKX_24}, the LDCU numerical flux was
developed for the Euler equations of gas dynamics, which, in the 1-D case, read as \eref{1.1} with $\mU=(\rho,\rho u,E)^\top$ and
$\mF=\big(\rho u,\rho u^2+p,u(E+p)\big)^\top$. Here, $\rho$, $u$, $p$, and $E$ are the density, velocity, pressure, and total energy,
respectively. The system is completed through the following equations of state (EOS) for ideal gases:
\begin{equation*}
p=(\gamma-1)\Big[E-\hf\rho u^2\Big],
\end{equation*}
where the parameter $\gamma$ represents the specific heat ratio.

The LDCU numerical fluxes are computed using the one-sided cell interface values $\mU^\pm_\jph$, which are reconstructed from the cell
averages $\big\{\,\xbar\mU_j\big\}$. In order to suppress oscillations, which may appear near shocks and contact discontinuities, we carry
out the reconstruction in local characteristic variables. To this end, we define $\widehat A_\jph:=A\big(\widehat\mU_\jph\big)$, and
compute the matrices $R_\jph$ and $R^{-1}_\jph$ such that $R^{-1}_\jph\widehat A_\jph R_\jph$ is diagonal. Here, $\widehat\mU_\jph$ is an
average of $\mU_j$ and $\mU_{j+1}$, which is obtained by averaging the primitive variables:
$$
\widehat\rho_\jph=\frac{\xbar\rho_j+\,\xbar\rho_{j+1}}{2},\quad\widehat u_\jph=\frac{u_j+u_{j+1}}{2},\quad
\widehat p_\jph=\frac{p_j+p_{j+1}}{2}.
$$
We then introduce the local characteristic variables in the neighborhood of $x=x_\jph$ (see, e.g., \cite{CCHKL_22,Joh,Qiu02,Shu20}):
\begin{equation*}
\bm\Gamma_\ell=R^{-1}_\jph\,\xbar\mU_{j+\ell},\quad\ell=-1,0,1,2.
\end{equation*}
Using $\bm\Gamma_{-1},\bm\Gamma_0,\bm\Gamma_1,\bm\Gamma_2$, we compute the slopes using the SBM-type limiters introduced in \cite{Lie03}:
\begin{equation}
(\bm\Gamma_x)_0=\phi^{\rm SBM}_{\theta,\tau}\!\left(\frac{\bm\Gamma_1-\bm\Gamma_0}{\bm\Gamma_0-\bm\Gamma_{-1}}\right)
\frac{\bm\Gamma_0-\bm\Gamma_{-1}}{\dx}\quad\mbox{and}\quad
(\bm\Gamma_x)_1=\phi^{\rm SBM}_{\theta,\tau}\!\left(\frac{\bm\Gamma_2-\bm\Gamma_1}{\bm\Gamma_1-\bm\Gamma_0}\right)
\frac{\bm\Gamma_1-\bm\Gamma_0}{\dx},
\label{2.10}
\end{equation}
where the function $\phi^{\rm SBM}_{\theta,\tau}$, defined by
\begin{equation}
\phi^{\rm SBM}_{\theta,\tau}(r):=\begin{cases}0,&\mbox{if }r<0,\\[0.2ex]
\min\big\{r\theta,\,1+\tau(r-1)\big\},&\mbox{if }0\le r\le1,\\[0.2ex]
r\,\phi^{\rm SBM}_{\theta,\tau}\!\left(\dfrac{1}{r}\right), &\mbox{otherwise},\end{cases}
\label{2.12}
\end{equation}
is applied in a componentwise manner. The parameters $\theta\in[1,2]$ and $\tau$ control the properties of the resulting reconstruction. In
general, larger values of $\theta$ lead to less dissipative, but typically more oscillatory, reconstructions. In all numerical examples in
\S\ref{sec4}, we take $\theta=2$, which can be safely used when the limiters are applied to the local characteristic variables. The
parameter $\tau$ determines whether the limiter is dissipative, compressive, or overcompressive: if $\tau\ge0.5$, then the limiter is
dissipative and contact discontinuities are typically smeared; if $0\le\tau<0.5$, then the limiter is compressive and contact waves are
resolved sharply, but smooth extrema may be slightly compressed, producing mild kinks; if $\tau<0$, then the limiter is overcompressive and
while contact waves typically remain sharp for long times, smooth profiles may develop artificial ${\cal O}(1)$ jumps and/or stair-like
structures. A particular choice made in \cite{CKM2025} was $\theta=2$, $\tau=0.5$ (which corresponds to the Minmod2 limiter) and $\theta=2$,
$\tau=-0.25$ in the ``smooth'' and ``rough'' areas, respectively. 

Once the slopes \eref{2.10} are available, we compute
$$
\bm\Gamma^-_\hf=\bm\Gamma_0+\frac{\dx}{2}(\bm\Gamma_x)_0,\quad\bm\Gamma^+_\hf=\bm\Gamma_1-\frac{\dx}{2}(\bm\Gamma_x)_1,
$$
and recover the interface values in conservative variables by $\mU^\pm_\jph=R_\jph\bm\Gamma^\pm_\hf$.
\begin{remark}
For the Euler equations of gas dynamics, explicit expressions and implementation details for $R_\jph$ and $R^{-1}_\jph$ can be found in
\cite[Appendix A]{CCHKL_22}.
\end{remark}

\subsection{One-Dimensional Adaptive Scheme}\label{sec22}
We now turn to the description of the proposed adaptive scheme. The key idea is to continuously vary the parameter $\tau\in [-0.25,0.5]$ in
\eref{2.12} according to the local smoothness of the computed solution.

To detect the local smoothness, we employ the SI from \cite{RL87} applied to the density. As in \cite{CKM2025}, we first define the
quantities
\begin{equation}
{\cal E}_j:=\frac{\big|\,\xbar\rho_{j+1}-2\,\xbar\rho_j+\,\xbar\rho_{j-1}\big|}
{\big|\,\xbar\rho_{j+1}-\,\xbar\rho_j\big|+\big|\,\xbar\rho_j-\,\xbar\rho_{j-1}\big|
+\varepsilon\big[\big|\,\xbar\rho_{j+1}\big|+2\big|\,\xbar\rho_j\big|+\big|\,\xbar\rho_{j-1}\big|\big]},
\label{2.12a}
\end{equation}
and then apply a local averaging to introduce
\begin{equation*}
\xbar{\cal E}_j:=\frac{1}{6}\big[{\cal E}_{j+1}+4{\cal E}_j+{\cal E}_{j-1}\big].
\end{equation*}
In \eref{2.12a}, the parameter $\varepsilon$ acts as a ``noise'' filter that prevents the refinement of small-amplitude wiggles caused by
loss of monotonicity. In all of the numerical experiments reported in \S\ref{sec41}, we have taken $\varepsilon=0.2$ as in \cite{CKM2025,RL87}.

Given $\xbar{\cal E}_j$, we prescribe $\tau_j$ as
\begin{equation*}
\tau_j=\begin{cases}\dfrac{1}{8}\Big[1+3\tanh\!\big(2000(\texttt{C}-\,\xbar{\cal E}_j)\big)\Big],&\mbox{if }\,\xbar{\cal E}_j<C,\\[1.5ex]
\dfrac{1}{8}\Big[1+3\tanh\!\big(300(\texttt{C}-\,\xbar{\cal E}_j)\big)\Big],&\mbox{otherwise},\end{cases}
\end{equation*}
where $\texttt{C}>0$ is a tunable parameter. The above formula yields $\tau_j\approx0.5$ in ``smooth'' regions, where
$\xbar{\cal E}_j\ll\texttt{C}$, $\tau_j\approx-0.25$ in ``rough'' regions, where $\xbar{\cal E}_j\gg\texttt{C}$, while providing a gradual
transition in a zone between ``smooth'' and ``rough'' areas. We then use $\tau=\tau_j$ for computing the SBM-type limiters in cell $I_j$.
\begin{remark}
The constant $\texttt{C}$ is a positive parameter to be selected for each problem at hand. The robustness of the shock detection strategy
depends on how sensitive the algorithm is to the choice of $\texttt{C}$. In our computations, $\texttt{C}$ is tuned on a coarse mesh, and
the same value is then used on finer meshes. This tuning procedure is demonstrated in Example 2 in \S\ref{sec41}.
\end{remark}
\begin{remark}
The adaptive strategy in \cite{CKM2025} employs a discontinuous switching between $-0.25$ and $0.5$, so that in \cite{CKM2025} we have used
\begin{equation*}
\tau_j=\begin{cases}-0.25,&\mbox{if }\,\xbar{\cal E}_j>\texttt{C},\\0.5,&\mbox{otherwise}.\end{cases} 
\end{equation*}
\end{remark}

\section{Two-Dimensional Scheme Adaptation Algorithm}\label{sec3}
In this section, we extend the 1-D adaptive strategy developed in \S\ref{sec2} to the 2-D hyperbolic systems of conservation laws
\eref{1.2}. We first briefly review the 2-D LDCU scheme from \cite{CKX_24} (\S\ref{sec31}) and then introduce the new 2-D adaptive scheme
(\S\ref{sec32}).

\subsection{2-D Low-Dissipation Central-Upwind (LDCU) Scheme}\label{sec31}
We introduce a uniform mesh consisting of the cells $I_{j,k}:=[x_\jmh,x_\jph]\times[y_\kmh,y_\kph]$, centered at $(x_j,y_k)$, where
$x_j=\big(x_\jmh+x_\jph\big)/2$, $y_k=\big(y_\kmh+y_\kph\big)/2$, $x_\jph-x_\jmh\equiv\dx$, and $y_\kph-y_\kmh\equiv\dy$. We denote by
$\,\xbar\mU_{j,k}$ the cell averages of $\mU$ over $I_{j,k}$,
\begin{equation*}
\xbar\mU_{j,k}\approx\frac{1}{\dx\,\dy}\iint\limits_{I_{j,k}}\mU(x,y,t)\,{\rm d}y\,{\rm d}x,
\end{equation*}
and assume that, at a given time level $t\ge0$, the values $\big\{\,\xbar\mU_{j,k}\big\}$ are available.

Following \cite{CKX_24}, the 2-D LDCU semi-discrete scheme evolves the cell averages by numerically solving the following system of ODEs:
\begin{equation}
\frac{{\rm d}\,\xbar\mU_{j,k}}{{\rm d}t}=-\frac{\bm{{\cal F}}_{\jph,k}-\bm{{\cal F}}_{\jmh,k}}{\dx}-
\frac{\bm{{\cal G}}_{j,\kph}-\bm{{\cal G}}_{j,\kmh}}{\dy},
\label{3.1}
\end{equation}
where $\bm{{\cal F}}_{\jph,k}=\bmF\big(\mU^-_{\jph,k},\mU^+_{\jph,k}\big)$ and
$\bm{{\cal G}}_{j,\kph}=\bmG\big(\mU^-_{j,\kph},\mU^+_{j,\kph}\big)$ are the LDCU numerical fluxes from \cite{CKX_24} developed for the 2-D
Euler equations of gas dynamics, which read as \eref{1.2} with $\mU=(\rho,\rho u,\rho v,E)^\top$,
$\mF=\big(\rho u,\rho u^2+p,\rho uv,u(E+p)\big)^\top$, and $\mG=\big(\rho v,\rho uv,\rho v^2+p,v(E+p)\big)^\top$. Here, $u$ and $v$ are the
$x$- and $y$-directional velocities, respectively, and the pressure is given by the ideal-gas EOS
$$
p=(\gamma-1)\Big[E-\frac{\rho}{2}(u^2+v^2)\Big].
$$
The LDCU numerical fluxes are computed using the one-sided point values $\mU^\pm_{\jph,k}$ and $\mU^\pm_{j,\kph}$, which can be
reconstructed in local characteristic variables via the local characteristic decomposition; see \cite[Appendix A]{CKM2025} for details.

\subsection{Two-Dimensional Adaptive Scheme}\label{sec32}
We now turn to the description of the proposed 2-D adaptive scheme. As in the 1-D case, the key idea is to smoothly vary the parameter
$\tau\in[-0.25,0.5]$ in the SBM-type limiters \eref{2.12} according to the local smoothness of the computed solution.

To detect the local smoothness, we employ the 2-D version of the SI from \cite{RL87} applied to the density. To this end, we first compute 
${\cal E}_{j,k}:=\sqrt{{\cal E}^1_{j,k}/{\cal E}^2_{j,k}}$, where
\begin{equation*}
\begin{aligned}
{\cal E}^1_{j,k}&:=\big(\,\xbar\rho_{j+1,k}-2\,\xbar\rho_{j,k}+\,\xbar\rho_{j-1,k}\big)^2+\big(\,\xbar\rho_{j,k+1}-2\,\xbar\rho_{j,k}+
\,\xbar\rho_{j,k-1}\big)^2,\\
{\cal E}^2_{j,k}&:=\Big(\big|\,\xbar\rho_{j+1,k}-\,\xbar\rho_{j,k}\big|+\big|\,\xbar\rho_{j,k}-\,\xbar\rho_{j-1,k}\big|
+\varepsilon\big[\big|\,\xbar\rho_{j+1,k}\big|+2\big|\,\xbar\rho_{j,k}\big|+\big|\,\xbar\rho_{j-1,k}\big|\big]\Big)^2\\
&\hspace*{0.2cm}+\Big(\big|\,\xbar\rho_{j,k+1}-\,\xbar\rho_{j,k}\big|+\big|\,\xbar\rho_{j,k}-\,\xbar\rho_{j,k-1}\big|
+\varepsilon\big[\big|\,\xbar\rho_{j,k+1}\big|+2\big|\,\xbar\rho_{j,k}\big|+\big|\,\xbar\rho_{j,k-1}\big|\big]\Big)^2,
\end{aligned}
\end{equation*}
with $\varepsilon$ playing the role of a noise filter; see \S\ref{sec22} for further discussion. We then apply a local averaging and
introduce
\begin{equation*}
\begin{aligned}
\xbar{\cal E}_{j,k}&:=\frac{1}{36}\Big[{\cal E}_{j-1,k-1}+{\cal E}_{j-1,k+1}+{\cal E}_{j+1,k-1}+{\cal E}_{j+1,k+1}\\
&\hspace*{1.1cm}+4\big({\cal E}_{j-1,k}+{\cal E}_{j,k-1}+{\cal E}_{j,k+1}+{\cal E}_{j+1,k}\big)+16{\cal E}_{j,k}\Big].
\end{aligned}
\end{equation*}
Finally, we select $\tau_{j,k}$ by the same smooth mapping as in the 1-D case,
\begin{equation*}
\tau_{j,k}=\begin{cases}\dfrac{1}{8}\Big[1+3\tanh\!\big(2000(\texttt{C}-\,\xbar{\cal E}_{j,k})\big)\Big],&
\mbox{if }\,\xbar{\cal E}_{j,k}<C,\\[1.5ex]
\dfrac{1}{8}\Big[1+3\tanh\!\big(300(\texttt{C}-\,\xbar{\cal E}_{j,k})\big)\Big],&\mbox{otherwise},\end{cases}
\end{equation*}
where $\texttt{C}>0$ is a tunable threshold parameter, and use $\tau=\tau_{j,k}$ in cell $I_{j,k}$ to evaluate the SBM-type limiters there. 

\section{Numerical Examples}\label{sec4}
In this section, we test the developed adaptive schemes on a series of numerical examples and compare their performance with that of the
adaptive LDCU schemes from \cite{CKM2025}. To this end, we apply the proposed adaptive schemes to the 1-D and 2-D Euler equations of gas
dynamics. For the sake of brevity, the adaptive LDCU schemes from \cite{CKM2025} and the proposed new adaptive schemes will be referred to
as the OLD and NEW schemes, respectively.

In all of the numerical experiments, the semi-discrete ODE systems \eref{2.1} and \eref{3.1} are integrated in time by the three-stage
third-order strong stability preserving (SSP) Runge-Kutta method (see, e.g., \cite{Gottlieb11,Gottlieb12}) with the CFL number set to $0.4$.
We take $\gamma=1.4$ in Examples 1--7 and $\gamma=5/3$ in Example 8.

\subsection{One-Dimensional Examples}\label{sec41}
{\bf Example 1---Shock-Density Wave Interaction Problem.} In the first example taken from \cite{SO89}, we consider the shock-density wave
interaction problem with the following initial data,
\begin{equation*}
(\rho,u,p)\Big|_{(x,0)}=\begin{cases}\bigg(\dfrac{27}{7},\dfrac{4\sqrt{35}}{9},\dfrac{31}{3}\bigg),&x<-4,\\[0.8ex]
(1+0.2\sin(5x),0,1),&x>-4,\end{cases}
\end{equation*}
which are prescribed in the computational domain $[-5,15]$ subject to the free boundary conditions.

\begin{figure}[ht!]
\centerline{\includegraphics[trim=0.7cm 0.2cm 0.9cm 0.6cm, clip, width=4.8cm]{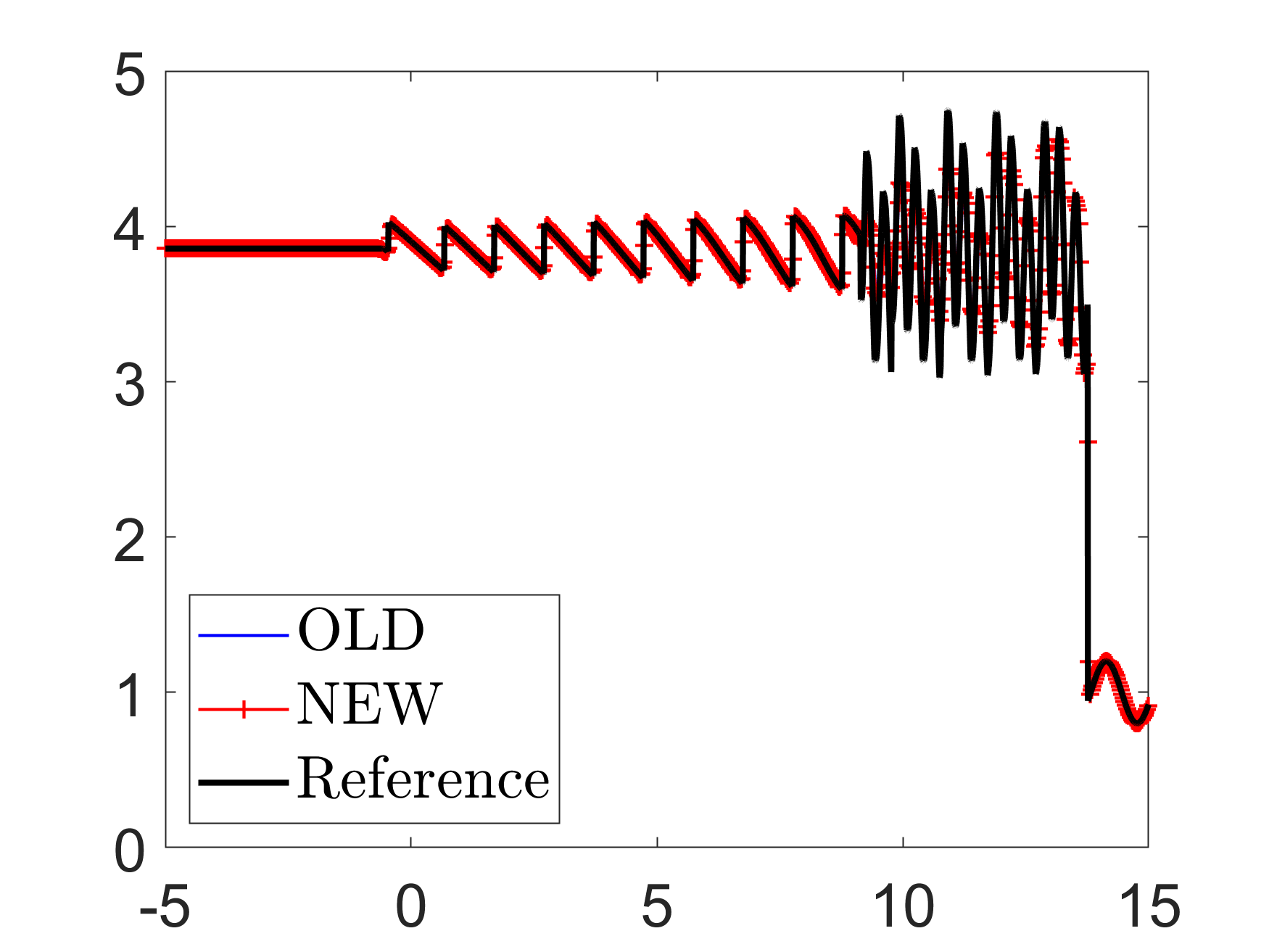}\hspace{0.8cm}
            \includegraphics[trim=0.7cm 0.2cm 0.9cm 0.6cm, clip, width=4.8cm]{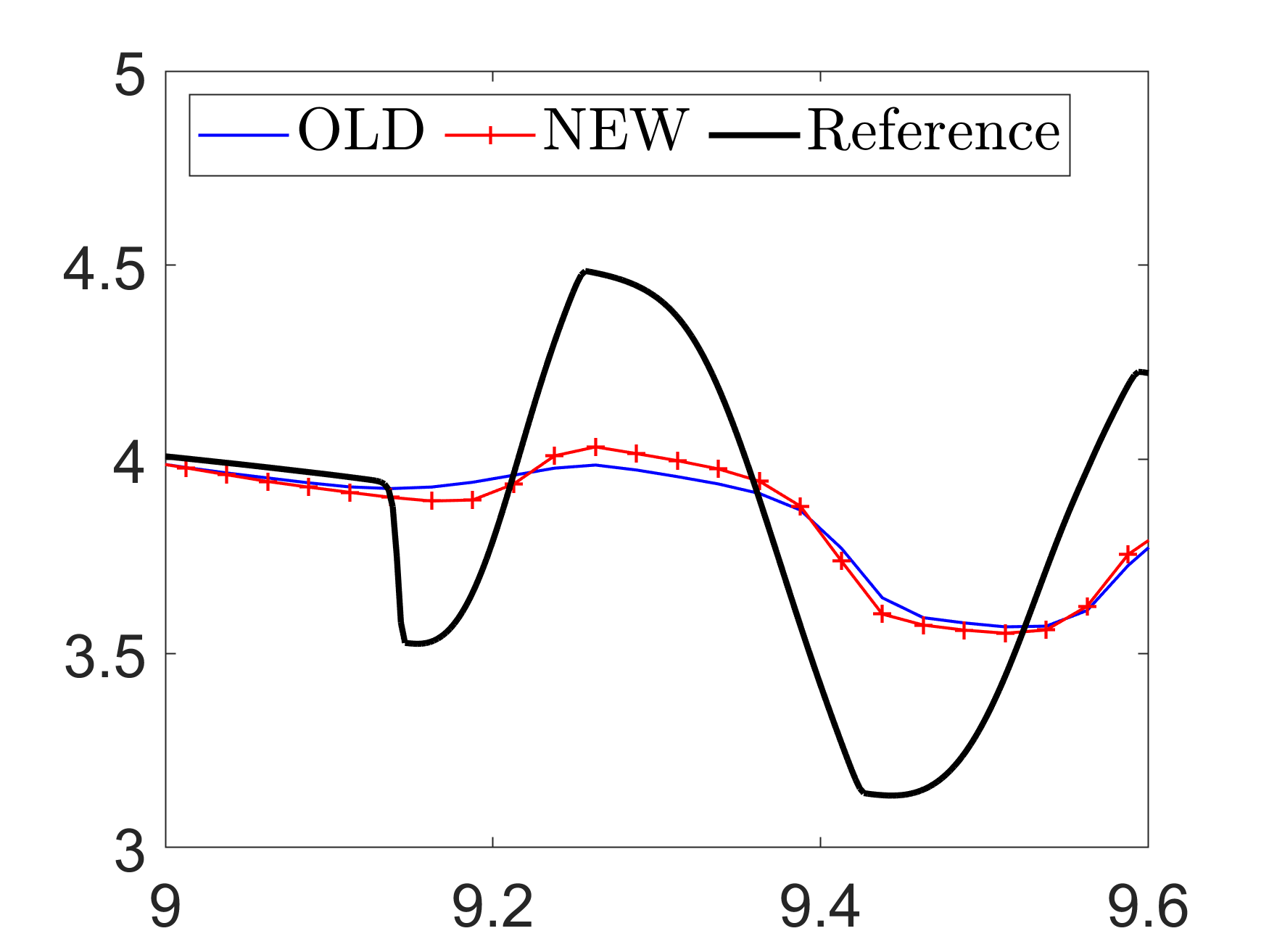}}
\caption{\sf Example 1: Density $\rho$ computed by the OLD and NEW schemes (left) and zoom at $x\in[9,9.6]$ (right).\label{fig4}}
\end{figure}

We compute the numerical solutions by the OLD (with the adaption constant $\texttt{C}=0.01$) and NEW (with the adaption constant
$\texttt{C}=0.005$) schemes on a uniform mesh with $\dx=1/40$ until the final time $t=5$, and plot the obtained numerical results in Figure
\ref{fig4} together with the reference solution computed by the LDCU scheme on a much finer mesh with $\dx=1/400$. It can be clearly seen
that the NEW scheme produces a slightly more accurate result compared to that obtained by the OLD scheme.

\noindent{\bf Example 2---Titarev-Toro Problem.} In the second example taken from \cite{Toro2005} (see also \cite{Shu88,Toro2005a}), we
consider the shock-entropy wave interaction problem. The initial conditions,
\begin{equation*}
(\rho,u,p)\Big|_{(x,0)}=\begin{cases}(1.51695,0.523346,1.805),&x<-4.5,\\(1+0.1\sin(20x),0,1),&x>-4.5,\end{cases}
\end{equation*}
correspond to a forward-facing shock wave of Mach number $1.1$ interacting with high-frequency density perturbations, that is, as the shock
wave moves, the perturbations spread ahead. The free boundary conditions are imposed in the computational domain $[-5,5]$.

\begin{figure}[ht!]
\centerline{\includegraphics[trim=0.7cm 0.2cm 1.2cm 0.6cm, clip, width=4.8cm]{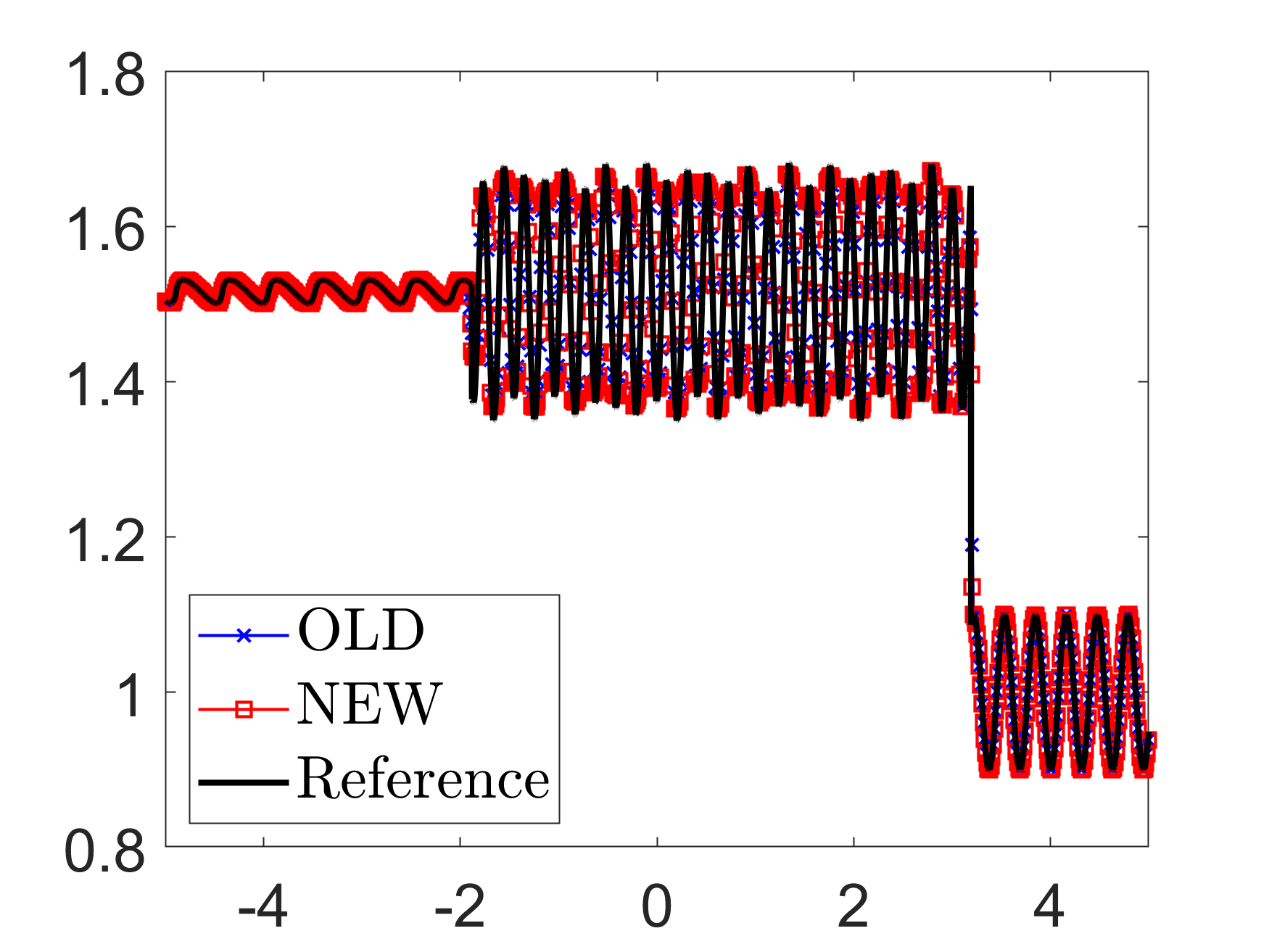}\hspace{0.8cm}
            \includegraphics[trim=0.7cm 0.2cm 1.2cm 0.6cm, clip, width=4.8cm]{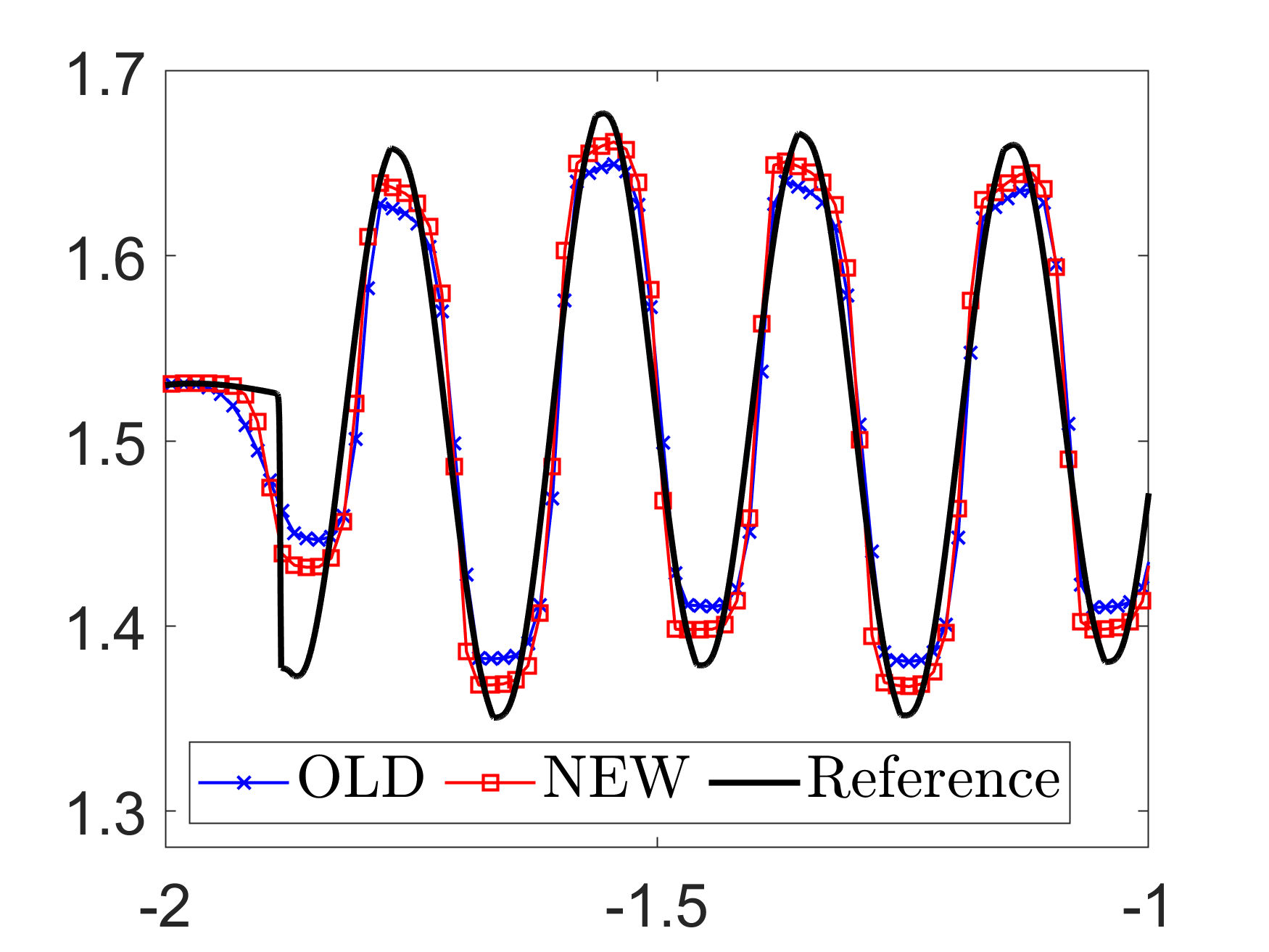}}
\caption{\sf Example 2: Density $\rho$ computed by the OLD and NEW schemes (left) and zoom at $x\in[-2,-1]$ (right).
\label{fig3}}
\end{figure}

We compute the numerical solution until the final time $t=5$ by the OLD (with the adaption constant $\texttt{C}=0.01$) and NEW (with the
adaption constant $\texttt{C}=0.002$) schemes on a uniform mesh with $\dx=1/80$. The obtained numerical results are plotted in Figure
\ref{fig3} along with the reference solution computed by the LDCU scheme on a much finer mesh with $\dx=1/1600$. In Figure \ref{fig3}
(right), we zoom the obtained solutions at the interval $[-2,-1]$, where one can see that this part of the solution is resolved much more
accurately by the NEW scheme.

\begin{figure}[ht!]
\centerline{\includegraphics[trim=0.7cm 0.2cm 1.4cm 0.3cm, clip, width=4.8cm]{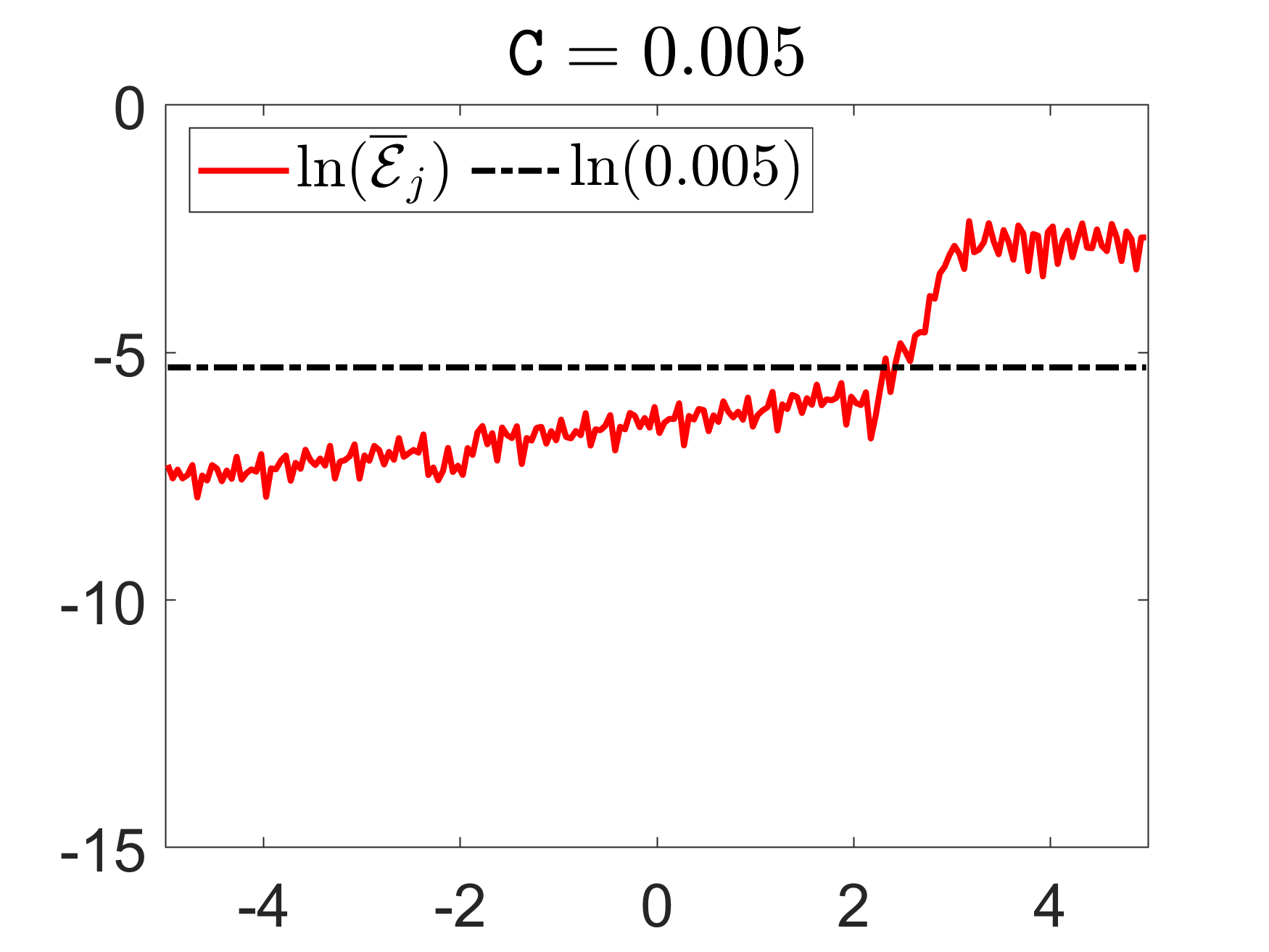}\hspace*{0.8cm}
            \includegraphics[trim=0.7cm 0.2cm 1.4cm 0.3cm, clip, width=4.8cm]{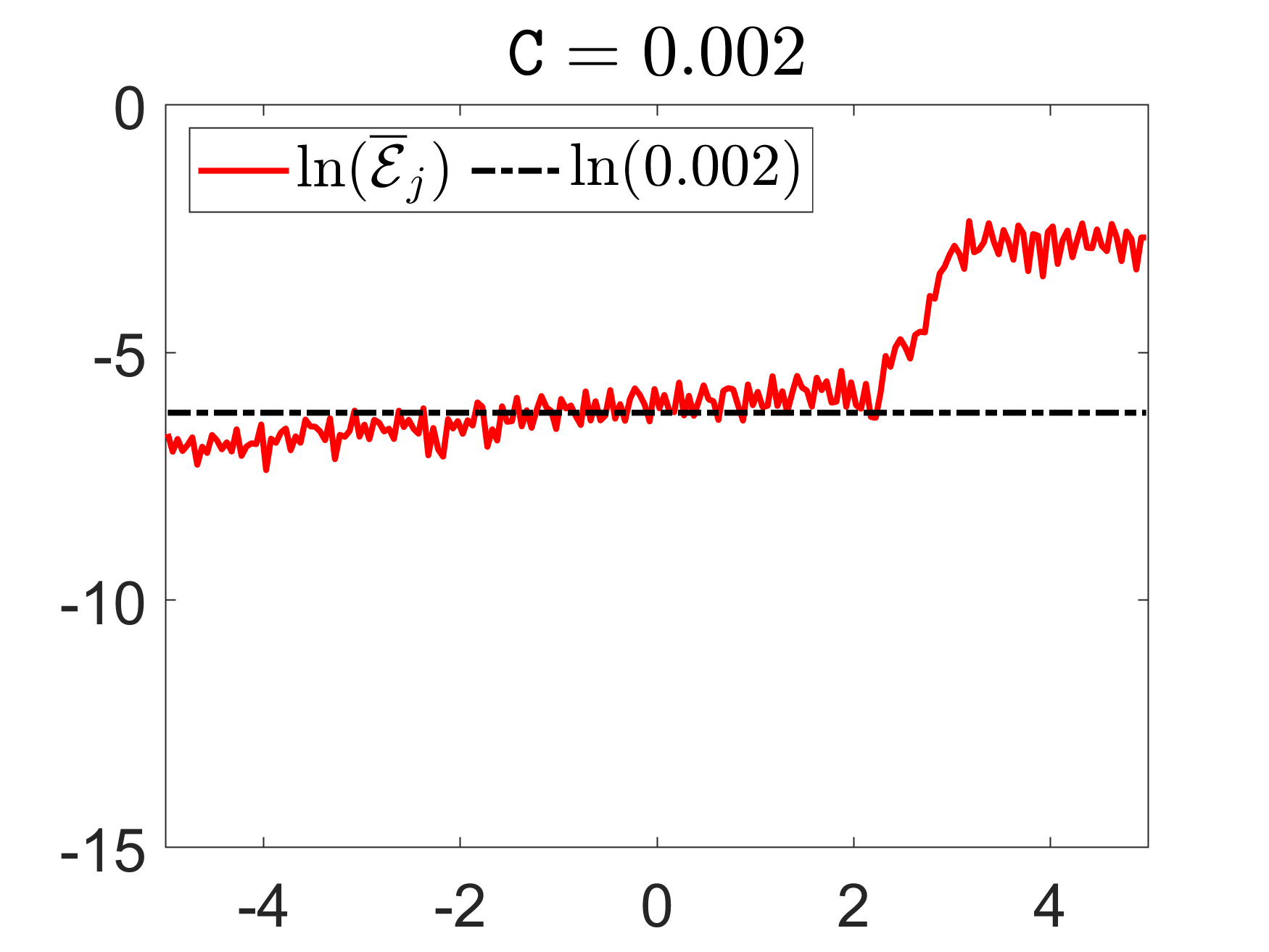}}
\vskip8pt
\centerline{\includegraphics[trim=0.7cm 0.2cm 1.4cm 0.3cm, clip, width=4.8cm]{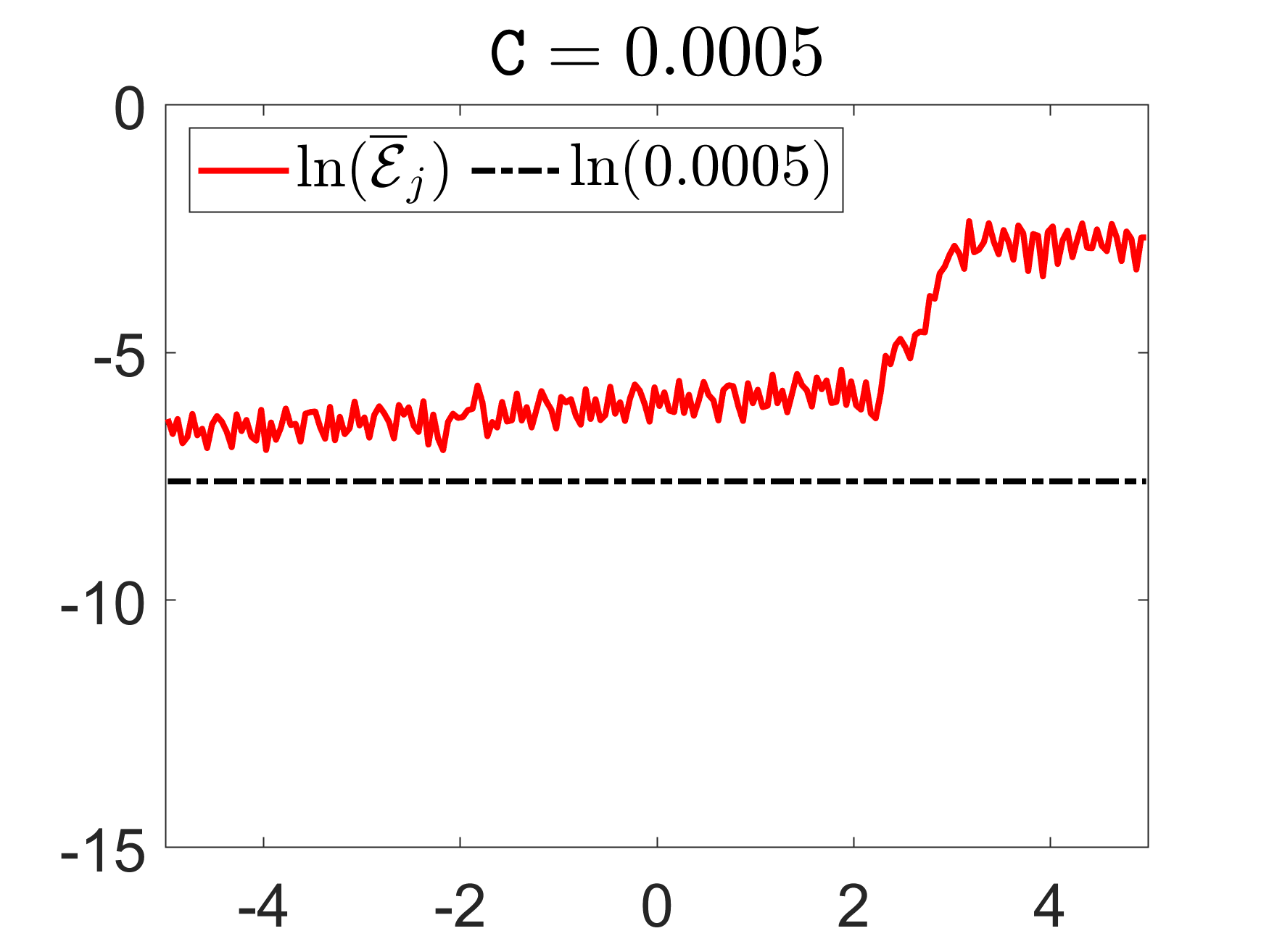}\hspace*{0.8cm}
     	    \includegraphics[trim=0.7cm 0.2cm 1.4cm 0.3cm, clip, width=4.8cm]{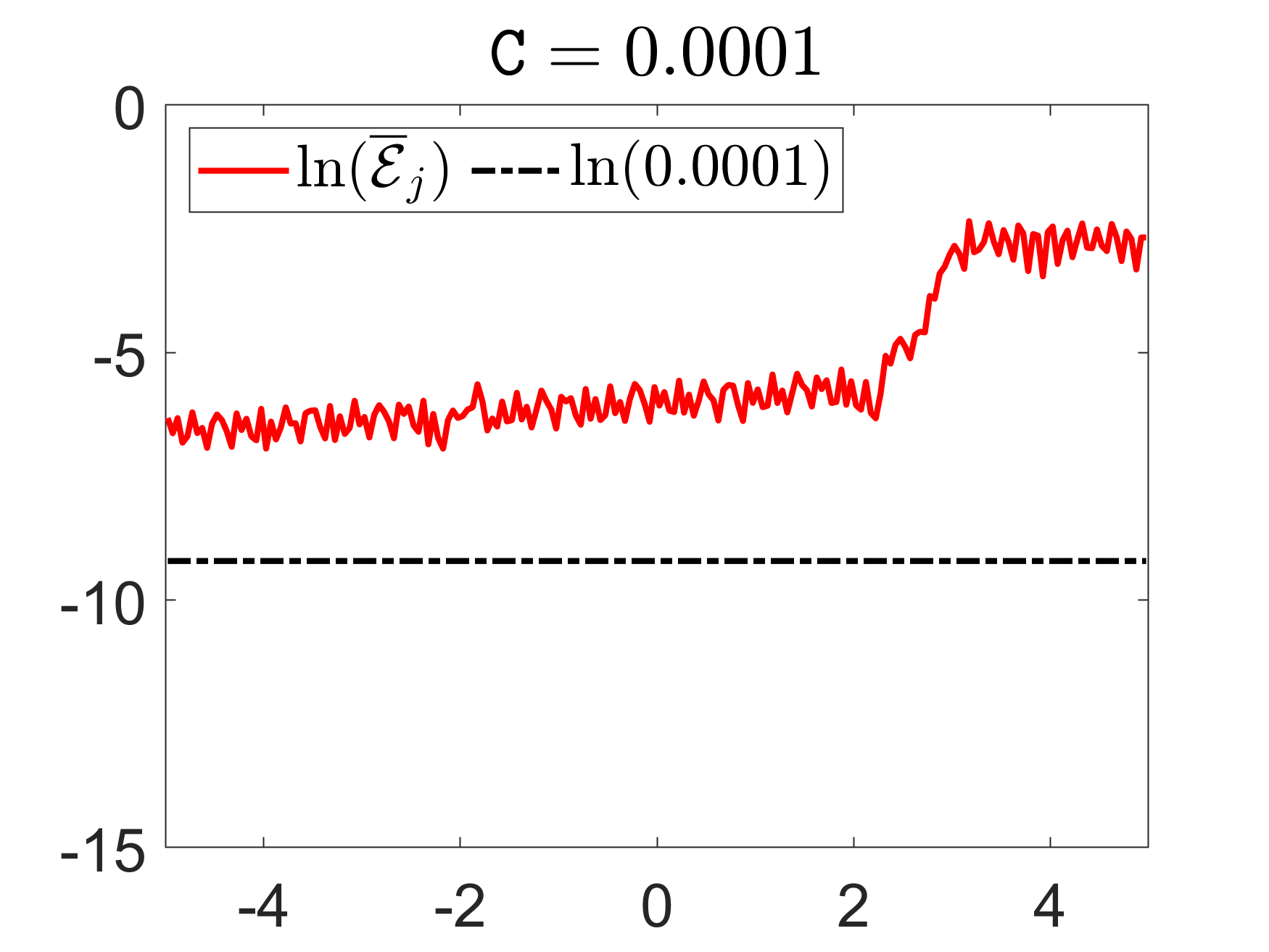}}
\caption{\sf Example 2: $\ln(\xbar{\cal E}_j)$ computed on a coarse mesh with $\dx=1/20$ and $\texttt{C}=0.005$ (top left),
$\texttt{C}=0.002$ (top right), $\texttt{C}=0.0005$ (bottom left), and $\texttt{C}=0.0001$ (bottom right).\label{fig4.33}}
\end{figure}

\begin{figure}[ht!]
\centerline{\includegraphics[trim=0.8cm 0.2cm 1.2cm 0.3cm, clip, width=4.8cm]{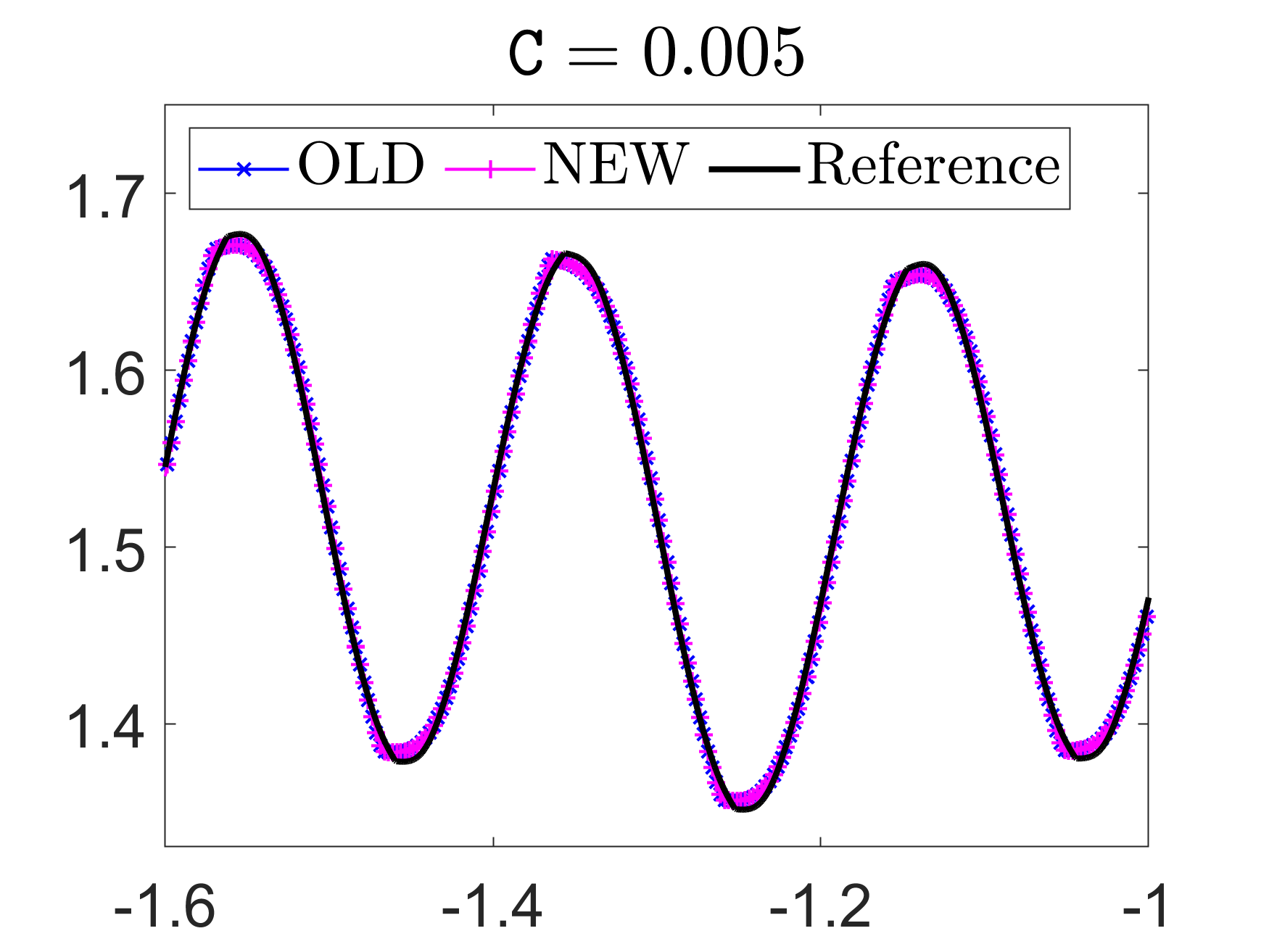}\hspace*{0.8cm}
            \includegraphics[trim=0.8cm 0.2cm 1.2cm 0.3cm, clip, width=4.8cm]{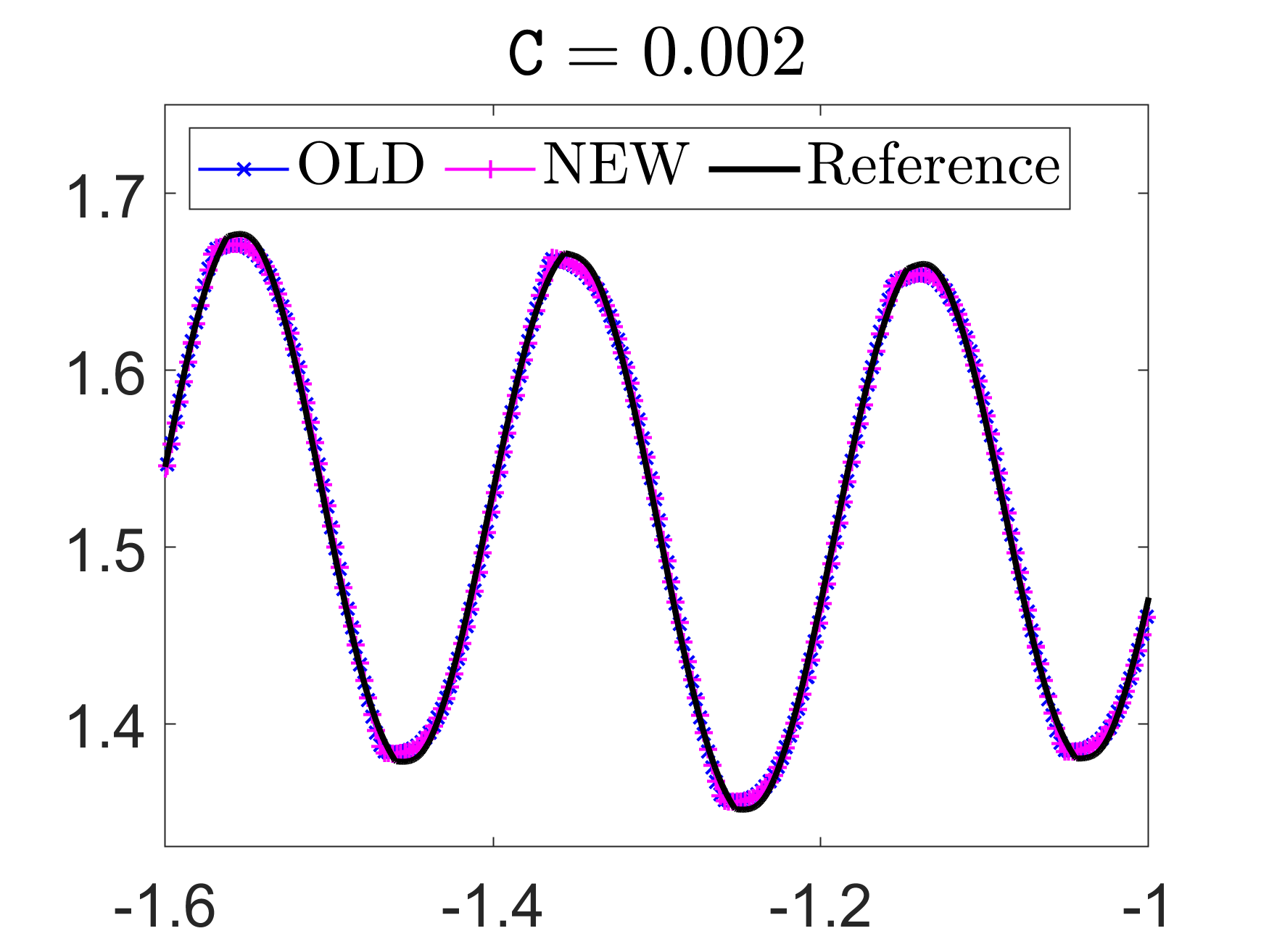}}
\vskip7pt
\centerline{\includegraphics[trim=0.8cm 0.2cm 1.2cm 0.3cm, clip, width=4.8cm]{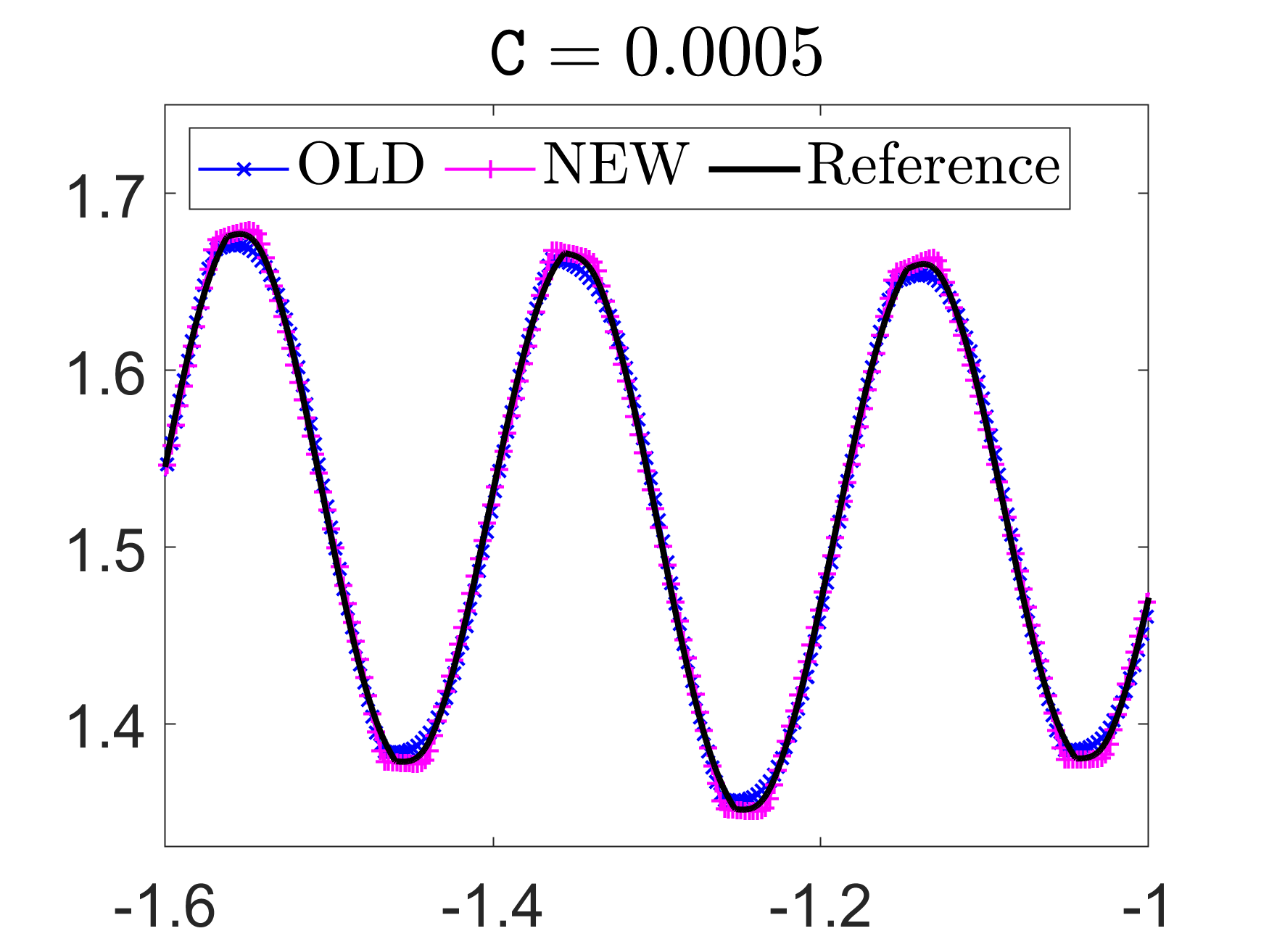}\hspace*{0.8cm}
     	    \includegraphics[trim=0.8cm 0.2cm 1.2cm 0.3cm, clip, width=4.8cm]{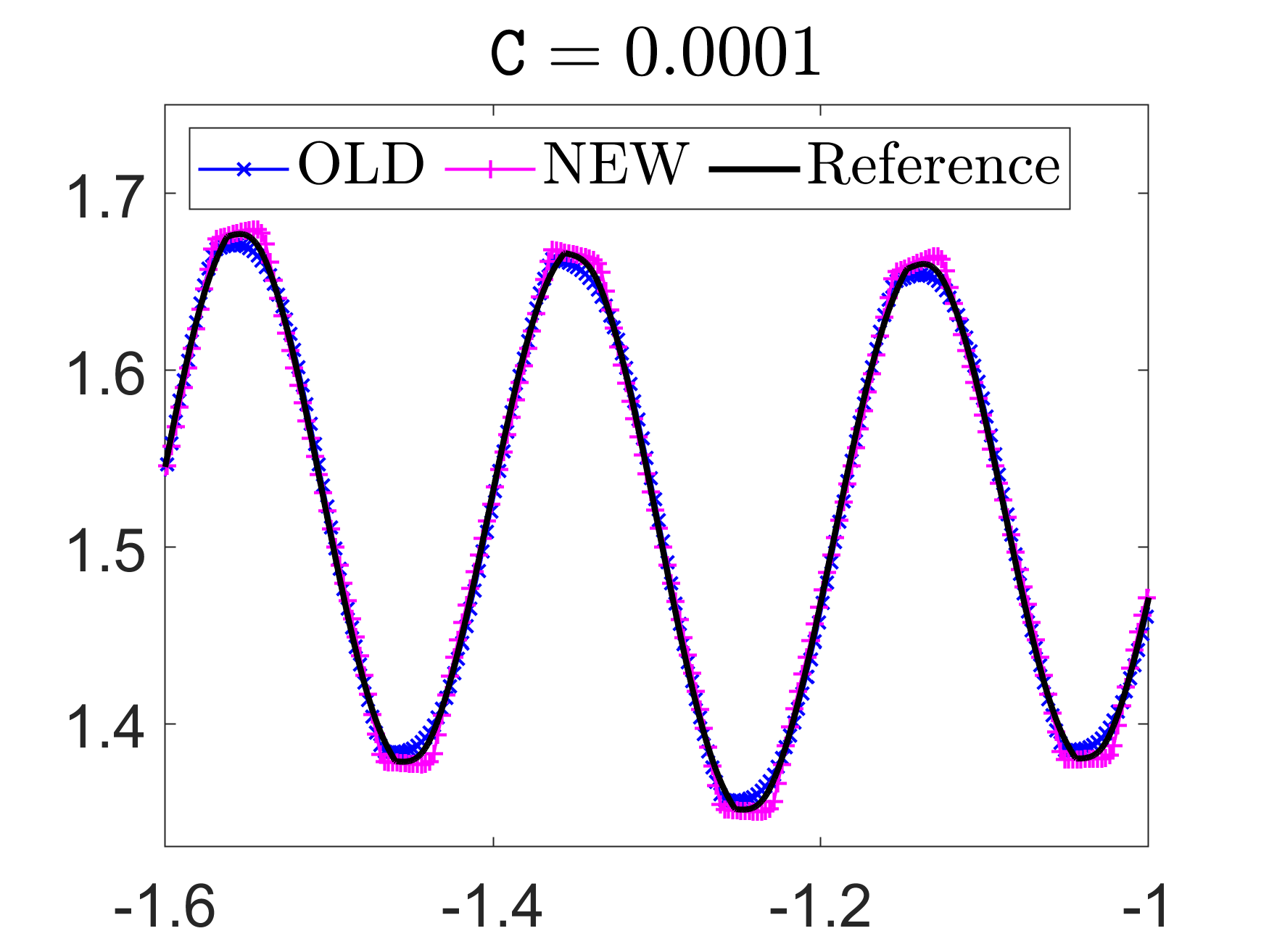}}
\caption{\sf Example 2: Density $\rho$ (zoom at $x\in[-1.6,-1]$) computed with $\dx=1/400$ and $\texttt{C}=0.005$ (top left),
$\texttt{C}=0.002$ (top right), $\texttt{C}=0.0005$ (bottom left), and $\texttt{C}=0.0001$ (bottom right).\label{fig4.44}}
\end{figure}

As mentioned above, the adaptive coefficient $\texttt{C}$ may be tuned as follows. First, it is adjusted on a coarse grid (this makes the
tuning process computationally inexpensive) and then used for high-resolution computation on finer grids. We demonstrate the tuning process
in this example. To this end, we first compute the numerical results on a coarse mesh with $\dx=1/20$ with $\texttt{C}=0.005$, $0.002$,
$0.0005$, and $0.0001$, and plot the corresponding values of $\ln({\,\xbar{\cal E}_j})$ computed at the final time step in Figure
\ref{fig4.33}. As one can see, the values $\texttt{C}=0.0005$ and $0.0001$ seem to be too small as for them the computed solution will be
identified everywhere as ``rough'' and the overcompressive limiter will be used throughout the entire computational domain. On the other
hand, the values $\texttt{C}=0.005$ and $0.002$ seem to be reasonable and the coarse mesh selection in Figure \ref{fig3} was thus produced
with $\texttt{C}=0.002$. We then compute the corresponding numerical results on a finer mesh with $\dx=1/400$ and plot the obtained
numerical results in Figure \ref{fig4.44}, where one can clearly see that the selection made on a coarse mesh is robust.

\bigskip 
\noindent{\bf Example 3---Blast Wave Problem.} In the last 1-D example, we consider a strong shock interaction problem from
\cite{Woodward88} subject to the following initial conditions:
\begin{equation*}
(\rho,u,p)\Big|_{(x,0)}=\begin{cases}(1,0,1000),&x<0.1,\\(1,0,0.01),&0.1\le x\le0.9,\\(1,0,100),&x>0.9,\end{cases}
\end{equation*}
which is considered on the interval $[0,1]$ subject to the solid wall boundary conditions.

\begin{figure}[ht!]
\centerline{\includegraphics[trim=1.3cm 0.2cm 1.3cm 0.8cm, clip, width=4.8cm]{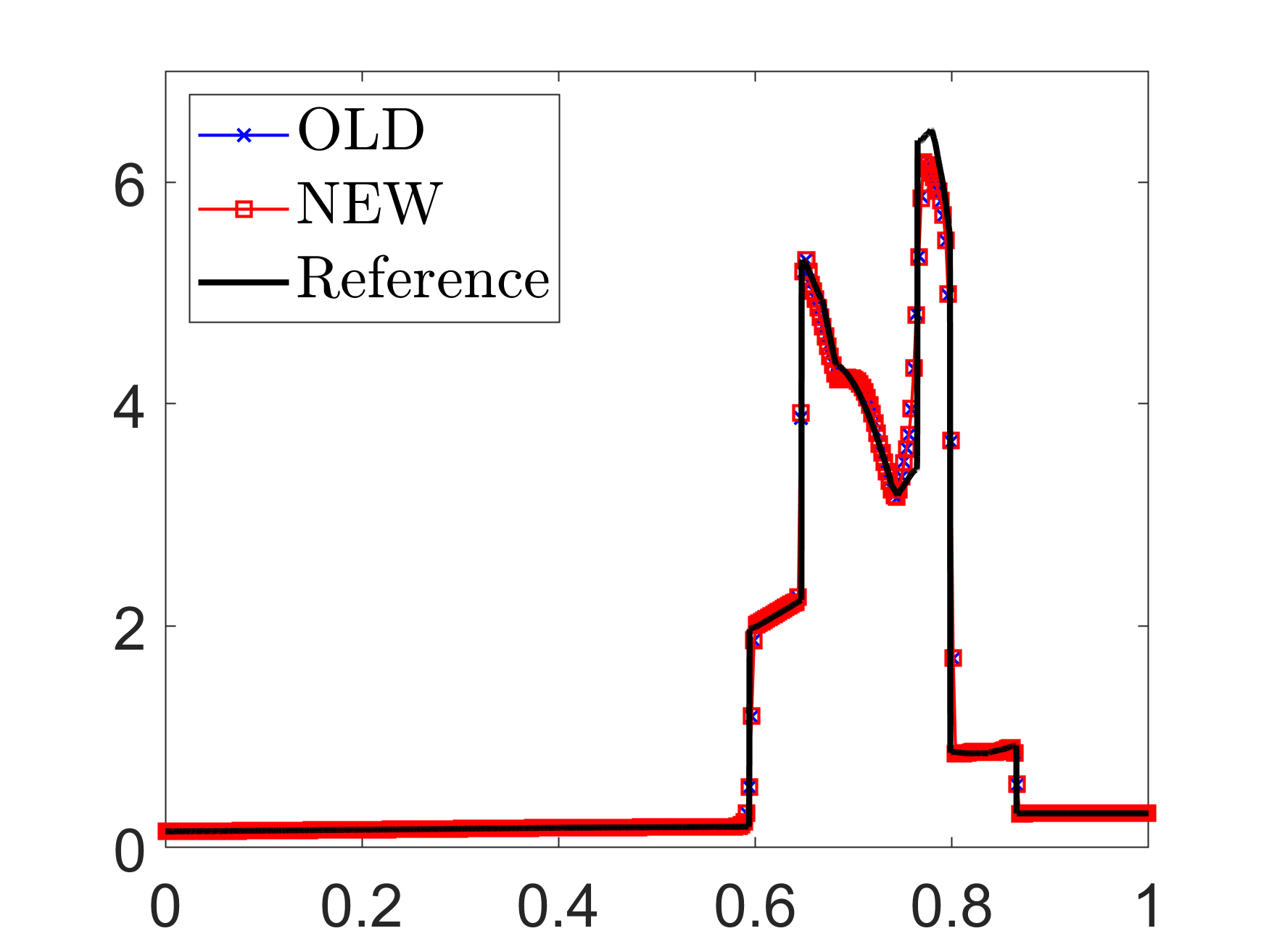}\hspace*{0.8cm}
            \includegraphics[trim=1.3cm 0.2cm 1.3cm 0.8cm, clip, width=4.8cm]{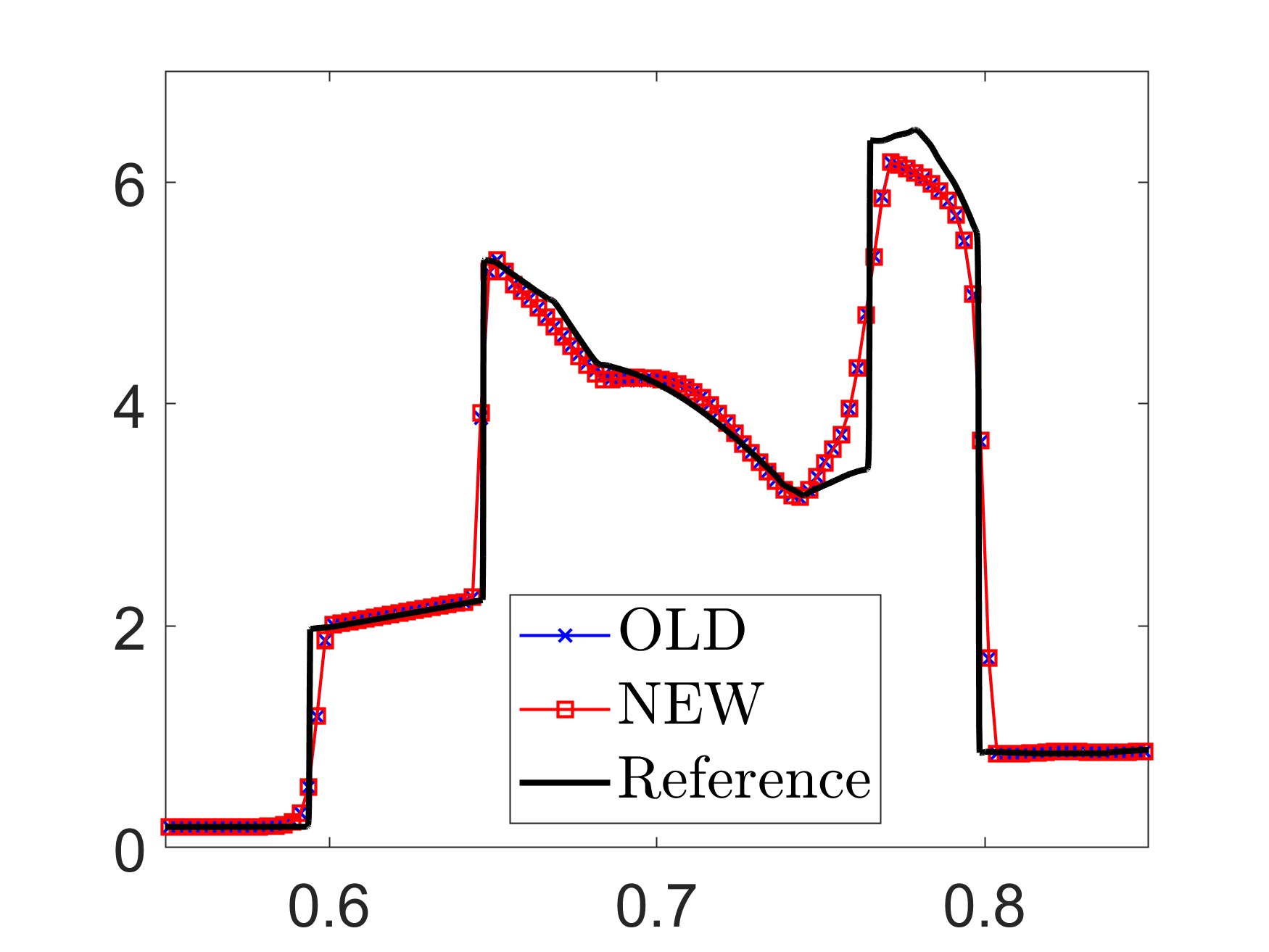}}
\caption{\sf Example 3: Density $\rho$ computed by the OLD and NEW schemes (left) and zoom at $x\in[0.55,0.85]$.\label{fig5}}
\end{figure}

We compute the numerical solutions until the final time $t=0.038$ by the OLD (with the adaption constant $\texttt{C}=0.01$) and NEW (with
the adaption constant $\texttt{C}=0.005$) schemes on a uniform mesh with $\dx=1/400$. The obtained results, together with the reference
solution computed by the LDCU scheme on a much finer mesh with $\dx=1/8000$, are plotted in Figure \ref{fig5}, where one can see that both
the OLD and NEW schemes produce similar results. This confirms the robustness of the new scheme adaption strategy: like the OLD scheme, the
NEW scheme is capable of achieving superb resolution of the contact wave located at about $x=0.6$. 

\subsection{Two-Dimensional Examples}
\noindent{\bf Example 4---2-D Accuracy Test.} In the first 2-D example taken from \cite{BD2013,Shu1998}, we consider the following smooth
initial data:
\begin{equation*}
\begin{aligned}
&\rho(x,y,0)=\bigg(1-\frac{(\gamma-1)\kappa^2}{2\gamma}\bigg)^{\frac{1}{\gamma-1}},\quad p(x,y,0)=\rho^\gamma(x,y,0),\\
	&u(x,y,0)=1-\kappa y, \quad v(x,y,0)=1+\kappa x, \quad \kappa = \frac{5}{2 \pi}e^{0.5(1-x^2-y^2)} 
\end{aligned}
\end{equation*}
subject to the periodic boundary conditions in the computational domain $[-10,10]\times[-10,10]$. The exact solution of this initial value
problem is given by $\mU(x,y,t)=\mU(x-t,y-t,0)$.

We compute the numerical solution until the final time $t=0.1$ by the NEW scheme (with the adaptation constant $\texttt{C}=0.05$) on a
sequence of uniform meshes with $\dx=\dy=1/10$, $1/20$, $1/40$, and $1/80$, and then measure the $L^1$-errors and the corresponding
experimental convergence rates for $\rho$, $u$, $v$, and $p$. The obtained results are presented in Table \ref{tab52}, where one can see
that the expected second order of accuracy has been achieved by the studied NEW scheme. 
\begin{table}[ht!]
\centering
\begin{tabular}{ccccccccc}
\toprule
\multirow{2}{*}{$\dx=\dy$}&\multicolumn{2}{c}{$\rho$}&\multicolumn{2}{c}{$u$}&\multicolumn{2}{c}{$v$}&\multicolumn{2}{c}{$p$}\\
\cline{2-9}
&Error&Rate&Error&Rate&Error&Rate&Error&Rate\\
\hline
$1/10$&3.67e-03&--- &6.07e-03&--- &6.18e-03&--- &4.46e-03&---\\
$1/20$&8.34e-04&2.14&1.47e-03&2.04&1.51e-03&2.04&1.00e-03&2.16\\
$1/40$&1.74e-04&2.26&3.22e-04&2.19&3.28e-04&2.20&2.07e-04&2.27\\
$1/80$&3.65e-05&2.25&6.74e-05&2.25&6.86e-05&2.26&4.28e-05&2.27\\
\bottomrule 
\end{tabular}
\caption{\sf Example 4: The $L^1$-errors of $\rho$, $u$, $v$, and $p$, and experimental convergence rates for the NEW scheme.
\label{tab52}}
\end{table}

\noindent{\bf Example 5---2-D Riemann Problem (Configuration 3).} In this example, we consider Configuration 3 of the 2-D Riemann
problems from \cite{Kurganov02} (also see \cite{Schulz93,Schulz93a,Zheng01}) with the initial conditions
\begin{equation*}
(\rho,u,v,p)\Big|_{(x,y,0)}=\begin{cases}(1.5,0,0,1.5),&x>1,~y>1,\\(0.5323,1.206,0,0.3),&x<1,~y>1,\\(0.138,1.206,1.206,0.029),&x<1,~y<1,\\
(0.5323,0,1.206,0.3),&x>1,~y<1,\end{cases}
\end{equation*}
prescribed in the computational domain $[0,1.2]\times[0,1.2]$ subject to the free boundary conditions.

\begin{figure}[ht!]
\centerline{\includegraphics[trim=0.cm 0cm 0cm 0cm, clip, width=10.cm]{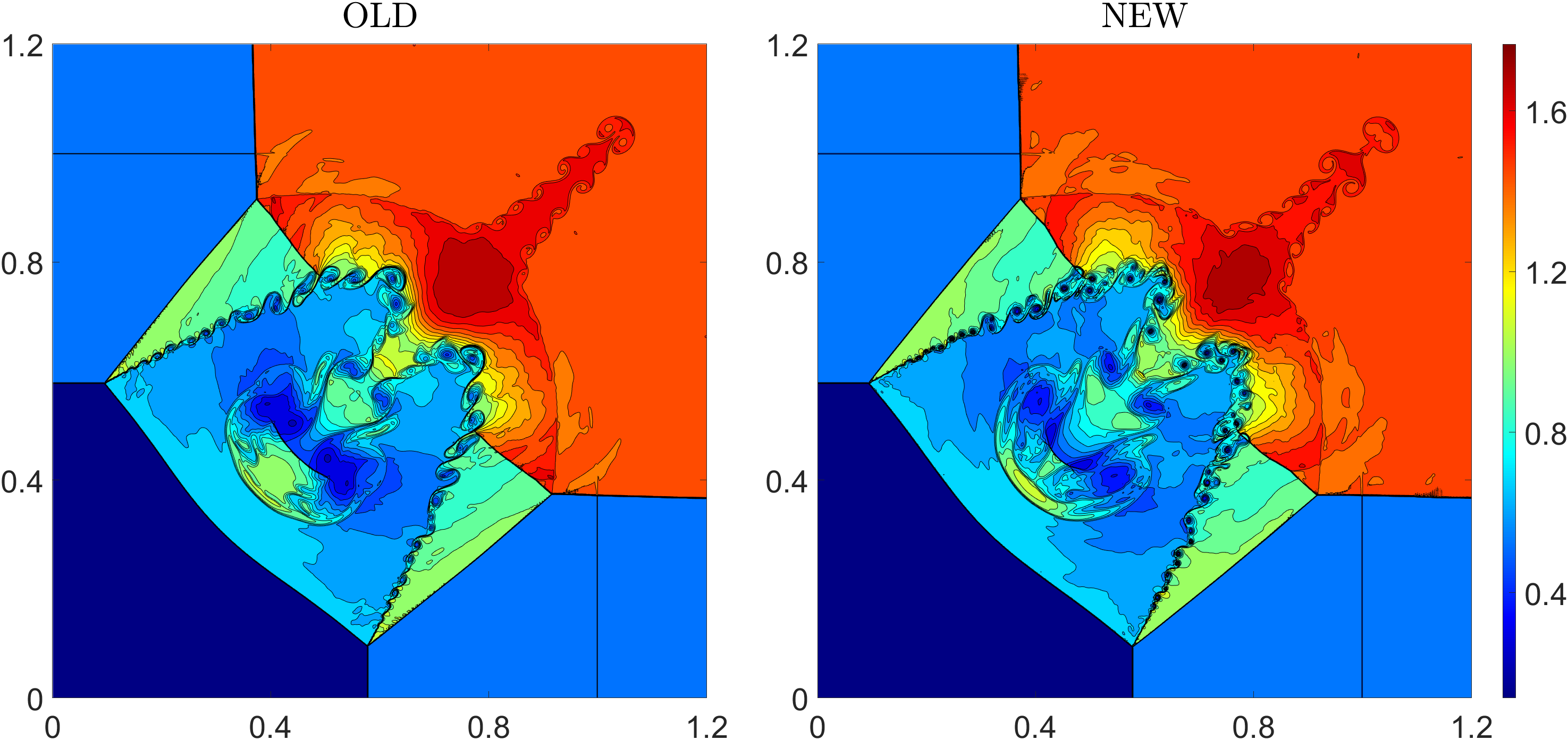}}
\caption{\sf Example 5: Density $\rho$ computed by the OLD (left) and NEW (right) schemes.\label{fig6a}}
\end{figure}

We compute the numerical solution until the final time $t=1$ by the OLD (with the adaption constant $\texttt{C}=0.08$) and NEW (with the
adaption constant $\texttt{C}=0.06$) schemes on a uniform mesh with $\dx=\dy=3/2500$ and present the obtained results in Figure \ref{fig6a},
where one can see that the NEW scheme slightly outperforms the OLD one in capturing a sideband instability of the jet in the zones of strong
along-jet velocity shear and the instability near the jet neck.

\noindent{\bf Example 6---2-D Riemann Problem (Configuration 6).} In this 2-D example, we consider Configuration 6 of the 2-D Riemann
problems from \cite{Kurganov02} with the initial conditions
\begin{equation*}
(\rho,u,v,p)\Big|_{(x,y,0)}=\begin{cases}(1,0.75,-0.5,1),&x>0.5,~y>0.5,\\(2,0.75,0.5,1),&x<0.5,~y>0.5,\\(1,-0.75,0.5,1),&x<0.5,~y<0.5,\\
(3,-0.75,-0.5,1),&x>0.5,~y<0.5,\end{cases}
\end{equation*}
which are prescribed in the computational domain $[0,1]\times[0,1]$ subject to the free boundary conditions.

\begin{figure}[ht!]
\centerline{\includegraphics[trim=0.cm 0cm 0cm 0cm, clip, width=10.cm]{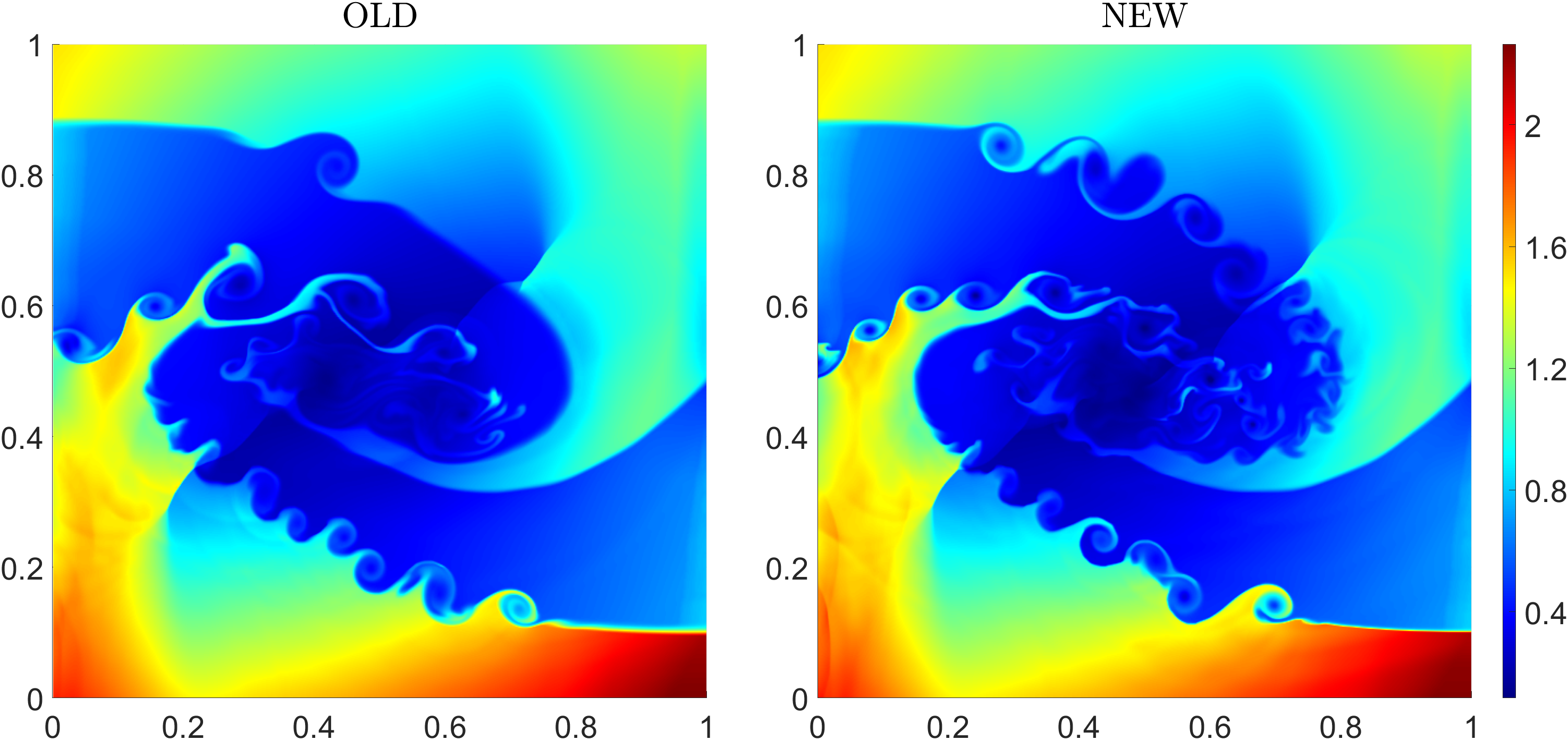}}
\caption{\sf Example 6: Density $\rho$ computed by the OLD (left) and NEW (right) schemes.\label{fig16a}}
\end{figure}

We compute the numerical solution until the final time $t=1$ by the OLD (with the adaption constant $\texttt{C}=0.1$) and NEW (with the
adaption constant $\texttt{C}=0.075$) schemes on a uniform mesh with $\dx=\dy=1/600$, and plot the obtained results in Figure \ref{fig16a}.
One can see that the NEW scheme is capable of capturing much more complicated vortex structures compared with the OLD one, which
demonstrates higher resolution and lower numerical dissipation of the NEW scheme.

\bigskip 
\noindent{\bf Example 7---2-D Riemann Problem (Configuration 12).} In this example, we consider Configuration 12 of the 2-D Riemann problems
from \cite{Kurganov02} with the initial conditions
\begin{equation*}
(\rho,u,v,p)\Big|_{(x,y,0)}=\begin{cases}(0.5313,0,0,0.4),&x>0.5,~y>0.5,\\(1,0.7276,0,1),&x<0.5,~y>0.5,\\(0.8,0,0,1),&x<0.5,~y<0.5,\\
(1,0,0.7276,1),&x>0.5,~y<0.5,\end{cases}
\end{equation*}
prescribed in the computational domain $[0,0.6]\times[0,0.6]$ subject to the free boundary conditions.

\begin{figure}[ht!]
\centerline{\includegraphics[trim=0.cm 0cm 0cm 0cm, clip, width=10.cm]{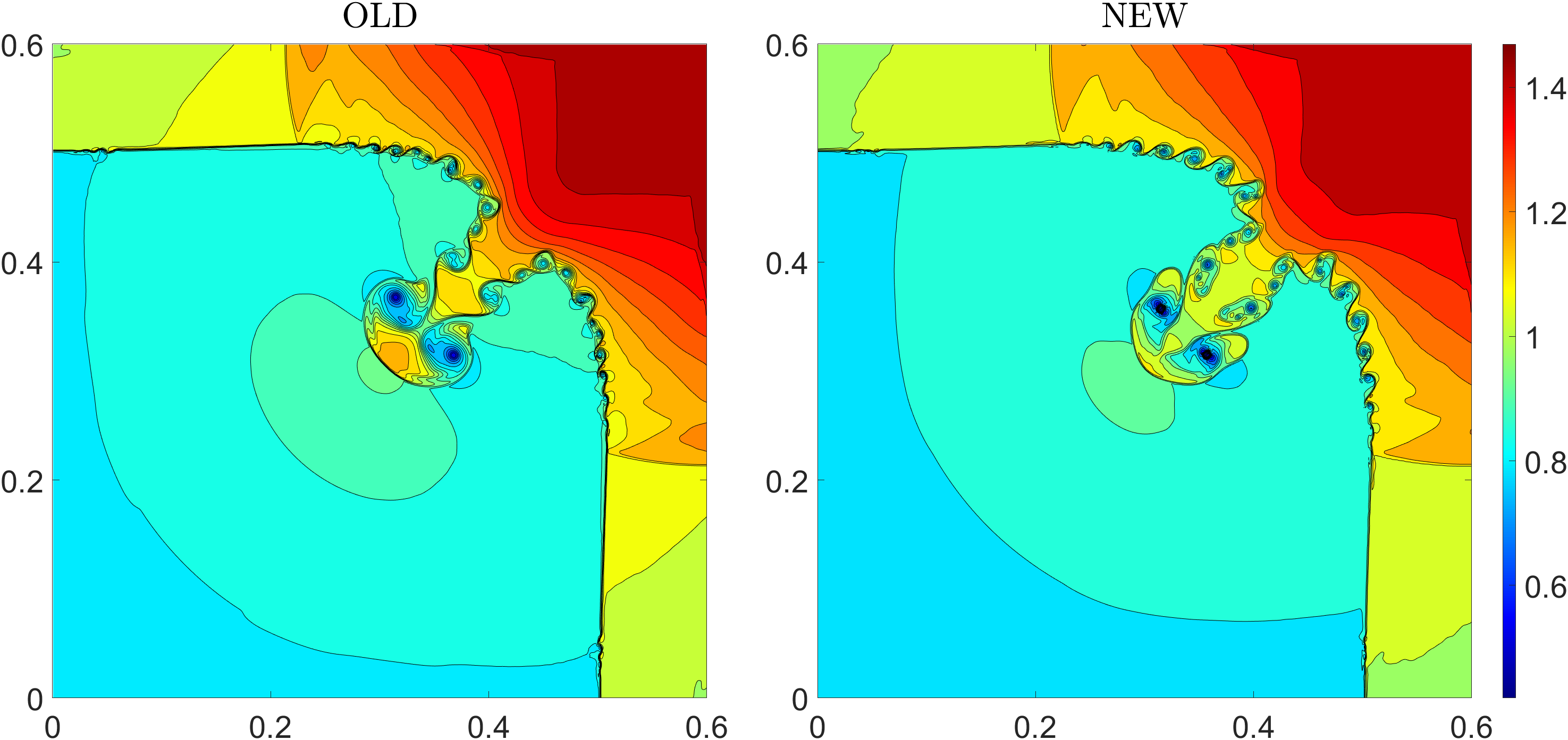}}
\caption{\sf Example 7: Density $\rho$ computed by the OLD (left) and NEW (right) schemes.\label{fig17a}}
\end{figure}

We compute the numerical solution until the final time $t=0.5$ by the OLD (with the adaption constant $\texttt{C}=0.03$) and NEW (with the
adaption constant $\texttt{C}=0.025$) schemes on a uniform mesh with $\dx=\dy=1/1000$, and plot the obtained results in Figure \ref{fig17a}.
As one can see, the NEW scheme resolves more vortices arising along the unstable contact surfaces compared with the OLD scheme. 

\noindent{\bf Example 8---RT Instability.} In the last example taken from \cite{Shi03}, we investigate the RT instability, which is a
physical phenomenon occurring when a layer of heavier fluid is placed on top of a layer of lighter fluid. To this end, we first modify the
2-D Euler equations of gas dynamics  by adding the gravitational source terms acting in the positive direction of the $y$-axis into the RHS
of the system:
\begin{equation*}
\begin{aligned}
&\rho_t+(\rho u)_x+(\rho v)_y=0,\\
&(\rho u)_t+(\rho u^2 +p)_x+(\rho uv)_y=0,\\
&(\rho v)_t+(\rho uv)_x+(\rho v^2+p)_y=\rho,\\
&E_t+\left[u(E+p)\right]_x+\left[v(E+p)\right]_y=\rho v,
\end{aligned}
\end{equation*}
and then use the following initial conditions:
\begin{equation*}
(\rho,u,v,p)\Big|_{(x,y,0)}=
\begin{cases}(2,0,-0.025c\cos(8\pi x),2y+1),&y<0.5,\\(1,0,-0.025c\cos(8\pi x),y+1.5),&\mbox{otherwise},\end{cases}
\end{equation*}
where $c:=\sqrt{\gamma p/\rho}$ is the speed of sound. The solid wall boundary conditions are imposed at $x=0$ and $x=0.25$, and the
following Dirichlet boundary conditions are specified at the top and bottom boundaries:
$$
(\rho,u,v,p)(x,1,t)=(1,0,0,2.5),\quad(\rho,u,v,p)(x,0,t)=(2,0,0,1).
$$

\begin{figure}[ht!]
\centerline{\includegraphics[trim=0.cm 0cm 0cm 0cm, clip, width=10.cm]{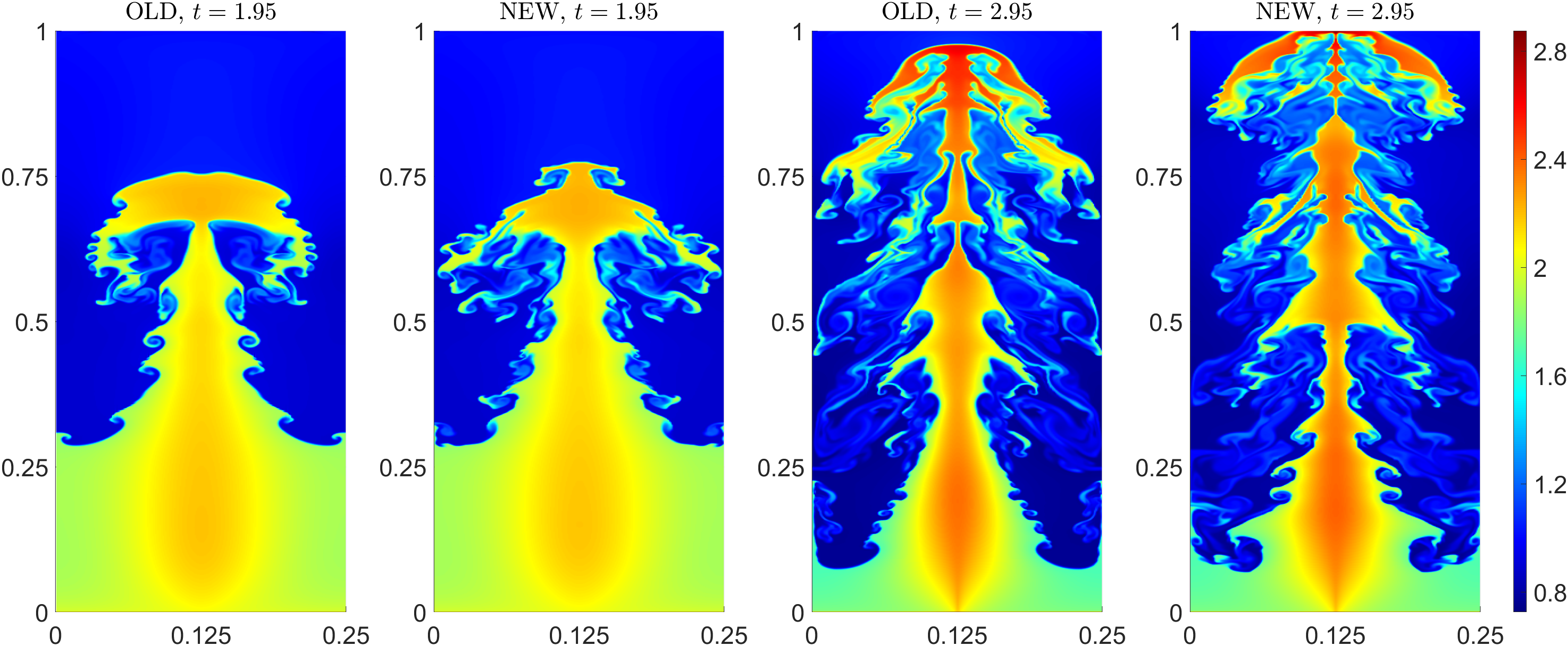}}
\caption{\sf Example 8: Density $\rho$ computed by the OLD and NEW schemes at $t=1.95$ and $2.95$.\label{fig10a}}
\end{figure}

We compute the numerical solution until the final time $t=2.95$ by the OLD (with the adaption constant $\texttt{C}=0.08$) and NEW (with the
adaption constant $\texttt{C}=0.06$) schemes on a uniform mesh with $\dx=\dy=1/1024$ in the computational domain $[0,0.25]\times[0,1]$. The
obtained numerical results are presented in Figure \ref{fig10a}, where one can see that the NEW scheme is less dissipative than the OLD one
in capturing the vortices.

\begin{acknowledgement}
The work of S. Chu and  M. Herty was funded by the Deutsche Forschungsgemeinschaft (DFG, German Research Foundation) - SPP 2410 Hyperbolic
Balance Laws in Fluid Mechanics: Complexity, Scales, Randomness (CoScaRa) within the Projects HE5386/26-1 (Numerische Verfahren für
gekoppelte Mehrskalenprobleme,525842915) and (Zufällige kompressible Euler Gleichungen: Numerik und ihre Analysis, 525853336) HE5386/27-1,
and  the Deutsche Forschungsgemeinschaft (DFG, German Research Foundation) - SPP 2183: Eigenschaftsgeregelte Umformprozesse with the
Projects HE5386/19-2,19-3 Entwicklung eines flexiblen isothermen Reckschmiedeprozesses für die eigenschaftsgeregelte Herstellung von
Turbinenschaufeln aus Hochtemperaturwerkstoffen (424334423). The work of A. Kurganov was supported in part by NSFC grant W2431004.
\end{acknowledgement}

\ethics{Competing Interests}{The authors have no conflicts of interest to declare that are relevant to the content of this chapter.}


\begin{thebibliography}{99.}%
%
%

\bibitem{BD2013} Boscheri, W., Dumbser, M.: Arbitrary-{L}agrangian-{E}ulerian one-step {WENO} finite volume schemes on unstructured triangular meshes. Commun. Comput. Phys. \textbf{14}(5), 1174--1206 (2013)

\bibitem{CCHKL_22} Chertock, A., Chu, S., Herty, M., Kurganov, A., Luk\'a\v{c}ov\'a-Medvi\v{d}ov\'a, M.: Local characteristic decomposition based central-upwind scheme. J. Comput. Phys. \textbf{473}, Paper No. 111718 (2023)

\bibitem{CKM2025} Chu, S., Kurganov, A., Menshov, I.: New adaptive low-dissipation central-upwind schemes. Appl. Numer. Math. \textbf{209}, 155--170 (2025)


\bibitem{CKX_24} Chu, S., Kurganov, A., Xin, R.: New low-dissipation central-upwind schemes. {P}art {II}. J. Sci. Comput. \textbf{103}(1), Paper No. 33 (2025)



\bibitem{Gottlieb11} Gottlieb, S., Ketcheson, D., Shu, C.-W.: Strong stability preserving {R}unge-{K}utta and multistep time discretizations. World Scientific Publishing Co. Pte. Ltd., Hackensack, NJ (2011)

\bibitem{Gottlieb12} Gottlieb, S., Shu, C.-W., Tadmor, E.: Strong stability-preserving high-order time discretization methods. SIAM Rev. \textbf{43}(1), 89--112 (2001)



\bibitem{Joh} Johnsen, E.: On the treatment of contact discontinuities using {WENO} schemes. J. Comput. Phys. \textbf{230}, 8665--8668 (2011)

\bibitem{Kurganov02} Kurganov, A., Tadmor, E.: Solution of two-dimensional {R}iemann problems for gas dynamics without {R}iemann problem solvers. Numer. Methods Partial Differential Equations \textbf{18}(5), 584--608 (2002)



\bibitem{Lie03} Lie, K.-A., Noelle, S.: On the artificial compression method for second-order nonoscillatory central difference schemes for systems of conservation laws. SIAM J. Sci. Comput. \textbf{24}(4), 1157--1174 (2003)

\bibitem{RL87} L\"ohner, R.: An adaptive finite element scheme for transient problems in {CFD}. Comput. Methods Appl. Mech. Eng. \textbf{61}(2), 323--338 (1987)

\bibitem{Qiu02} Qiu, J., Shu, C.-W.: On the construction, comparison, and local characteristic decomposition for high-order central {WENO} schemes. J. Comput. Phys. \textbf{183}(1), 187--209 (2002)

\bibitem{Schulz93} Schulz-Rinne, C. W.: Classification of the {R}iemann problem for two-dimensional gas dynamics. SIAM J. Math. Anal. \textbf{24}, 76--88 (1993)

\bibitem{Schulz93a} Schulz-Rinne, C. W., Collins, J. P., Glaz, H. M.: Numerical solution of the {R}iemann problem for two-dimensional gas dynamics. SIAM J. Sci. Comput. \textbf{14}, 1394--1414 (1993)

\bibitem{Shi03} Shi, J., Zhang, Y.-T., Shu, C.-W.: Resolution of high order {WENO} schemes for complicated flow structures. J. Comput. Phys. \textbf{186}, 690--696 (2003)

\bibitem{Shu20} Shu, C.-W.: Essentially non-oscillatory and weighted essentially non-oscillatory schemes. Acta Numer. \textbf{29}, 701--762 (2020)


\bibitem{Shu1998} Shu, C.-W.: Essentially non-oscillatory and weighted essentially non-oscillatory schemes for hyperbolic conservation laws. In: Advanced Numerical Approximation of Nonlinear Hyperbolic Equations (Cetraro, 1997), Lecture Notes in Math., vol. 1697, pp. 325--432. Springer, Berlin (1998)

\bibitem{Shu88} Shu, C.-W., Osher, S.: Efficient implementation of essentially non-oscillatory shock-capturing schemes. J. Comput. Phys. \textbf{77}, 439--471 (1988)

\bibitem{SO89} Shu, C.-W., Osher, S.: Efficient implementation of essentially nonoscillatory shock-capturing schemes. {II}. J. Comput. Phys. \textbf{83}(1), 32--78 (1989)

\bibitem{Toro2005} Titarev, V. A., Toro, E. F.: W{ENO} schemes based on upwind and centred {TVD} fluxes. Comput. \& Fluids. \textbf{34}, 705--720 (2005)


\bibitem{Toro2005a} Toro, E. F., Titarev, V. A.: T{VD} fluxes for the high-order {ADER} schemes. J. Sci. Comput. \textbf{24}, 285--309 (2005)

\bibitem{Woodward88} Woodward, P., Colella, P.: The numerical solution of two-dimensional fluid flow with strong shocks. J. Comput. Phys. \textbf{54}(1), 115--173 (1984)

\bibitem{Zheng01} Zheng, Y.: Systems of conservation laws. Two-dimensional Riemann problems. Birkh\"{a}user Boston, Inc., Boston, MA (2001)






  
    



\end{thebibliography}
\end{document}